\documentclass[12pt,english,a4paper]{amsart}
\usepackage[T1]{fontenc}
\usepackage[latin1]{inputenc}
\usepackage{geometry}
\usepackage{graphicx}
\usepackage{amssymb}

\makeatletter

\let\SF@@footnote\footnote
\def\footnote{\ifx\protect\@typeset@protect
    \expandafter\SF@@footnote
  \else
    \expandafter\SF@gobble@opt
  \fi
}
\expandafter\def\csname SF@gobble@opt \endcsname{\@ifnextchar[
  \SF@gobble@twobracket
  \@gobble
}
\edef\SF@gobble@opt{\noexpand\protect
  \expandafter\noexpand\csname SF@gobble@opt \endcsname}
\def\SF@gobble@twobracket[#1]#2{}

 \theoremstyle{plain}    
 \newtheorem{thm}{Theorem}[section]
 \numberwithin{equation}{section} 
 \numberwithin{figure}{section} 
 \theoremstyle{plain}
 \theoremstyle{plain}    
 \newtheorem{cor}[thm]{Corollary} 
 \theoremstyle{remark}
 \newtheorem{rem}[thm]{Remark}
 \theoremstyle{plain}    
 \newtheorem{conjecture}[thm]{Conjecture} 
 \theoremstyle{remark}    
 \newtheorem{acknowledgement}[thm]{Acknowledgement} 
 \theoremstyle{definition}
 \newtheorem{defn}[thm]{Definition}
 \theoremstyle{plain}    
 \newtheorem{lem}[thm]{Lemma} 
 \theoremstyle{plain}    
 \newtheorem{prop}[thm]{Proposition} 
 \theoremstyle{definition}
  \newtheorem{example}[thm]{Example}

\def\makebbb#1{
    \expandafter\gdef\csname#1\endcsname{
        \ensuremath{\Bbb{#1}}}
}
\makebbb{R}
\makebbb{N}
\makebbb{Z}
\makebbb{C}
\makebbb{G}
\makebbb{E}
\makebbb{H}
\makebbb{P}
\makebbb{B}
\makebbb{K}

\usepackage{babel}
\makeatother
\begin{document}

\title{Bergman kernels and equilibrium measures for line bundles over projective
manifolds }

\author{Robert Berman}

\email{robertb@math.chalmers.se}

\curraddr{Institut Fourier, 100 rue des Maths, BP 74, 38402 St Martin d'Heres
(France)}


\begin{abstract}
Let $L$ be a holomorphic line bundle over a compact complex projective
Hermitian manifold $X.$ Any fixed smooth hermitian metric $\phi$
on $L$ induces a Hilbert space structure on the space of global holomorphic
sections with values in the $k$th tensor power of $L.$ In this paper
various convergence results are obtained for the corresponding Bergman
kernels (i.e. orthogonal projection kernels). The convergence is studied
in the large $k$ limit and is expressed in terms of the equilibrium
metric $\phi_{e}$ associated to the fixed metric $\phi,$ as well
as in terms of the Monge-Ampere measure of the metric $\phi$ itself
on a certain support set. It is also shown that the equilibrium metric
is $\mathcal{C}^{1,1}$ on the complement of the augmented base locus
of $L.$ For $L$ ample these results give generalizations of well-known
results concerning the case when the curvature of $\phi$ is globally
positive (then $\phi_{e}=\phi$). In general, the results can be seen
as local metrized versions of Fujita's approximation theorem for the
volume of $L.$ 

\tableofcontents{}
\end{abstract}
\maketitle

\section{\label{sec:Introduction}Introduction}

Let $L$ be a holomorphic line bundle over a compact complex projective
manifold $X$ of dimension $n.$ Fix a smooth Hermitian fiber metric,
denoted by $\phi,$ on $L$ and a smooth volume form $\omega_{n}$
on $X.$ The curvature form of the metric $\phi$ may be written as
$dd^{c}\phi$ (see section \ref{sub:Notation-and-setup} for definitions
and further notation). Denote by $\mathcal{H}(X,L^{k})$ the Hilbert
space obtained by equipping the space $H^{0}(X,L^{k})$ of global
holomorphic sections with values in the tensor power $L^{k}$ with
the norm induced by the given metric $\phi$ on $L$ and the volume
form $\omega_{n}.$ The \emph{Bergman kernel} of the Hilbert space
$\mathcal{H}(X,L^{k})$ is the integral kernel of the orthogonal projection
from the space of all smooth sections with values in $L^{k}$ onto
$\mathcal{H}(X,L^{k}).$ It may be represented by a holomorphic section
$K_{k}(x,y)$ of the pulled back line bundle $L^{k}\boxtimes\overline{L}^{k}$
over $X\times\overline{X}$ (formula \ref{def: K}).

In the case when the curvature form $dd^{c}\phi$ is globally \emph{positive}
the asymptotic properties of the Bergman kernel $K_{k}(x,y)$ as $k$
tends to infinity have been studied thoroughly with numerous applications
in complex geometry and mathematical physics. For example, $K_{k}(x,y)$
admits a complete local asymptotic expansion in powers of $k;$ the
Tian-Zelditch-Catlin expansion (see \cite{ze,b-b-s} and references
therein). The point is that when the curvature form $dd^{c}\phi$
is globally positive, the Bergman kernel asymptotics at a fixed point
may be localized and hence only depend (up to negligable terms) on
the covariant derivatives of $dd^{c}\phi$ at the fixed point. 

The aim of the present paper is to study the case of a general smooth
metric $\phi$ on an arbitirary line bundle $L$ over a projective
manifold, where global effects become important and where there appears
to be very few previous general results even in the case when the
line bundle $L$ is ample (and even when $X$ is a complex curve).
We will be mainly concerned with three natural \emph{positive} measures
on $X$ associated to the setup introduced above. In order to introduce
these measures first assume that the line bundle $L$ is ample. The
first measure on $X$ to be considered is the \emph{equilibrium measure}\[
\mu_{\phi}:=(dd^{c}\phi_{e})^{n}/n!,\]
 where \emph{}$\phi_{e}$ is the \emph{equilibrium metric} defined
by the upper envelope \ref{eq:extem metric} (i.e. $\phi_{e}(x)=\sup\widetilde{\phi}(x),$
where the supremum is taken over all metrics $\widetilde{\phi}\leq\phi$
with positive curvature). For example, when $X$ is the projective
line $\P^{1}$ and $L$ is the hyperplane line bundle $\mathcal{O}(1)$
the measure $\mu_{\phi}$ is a minimizer of the {}``weighted logarithmic
energy'' \cite{s-t}. Next, the weak large $k$ limit of the measures
\begin{equation}
k^{-n}B_{k}\omega_{n},\label{eq:intro B meas}\end{equation}
 where $B_{k}(x):=K_{k}(x,x)e^{-k\phi}$ will be referred to as the
\emph{Bergman function} is considered and finally the limit of the
measure \[
(dd^{c}(k^{-1}\textrm{ln\,$K_{k}(x,x)))^{n}/n!$, }\]
often referred to as the $k$th \emph{Bergman volume form} on $X$
associated to $(L,\phi)$ (and $k^{-1}\textrm{ln\,$K_{k}(x,x)$ }$is
called the $k$th \emph{Bergman metric} on $L).$ 

When $L$ is ample it is well-known that the integrals over $X$ of
all three measures coincide. In fact, the integrals all equal the
integral over $X$ of the possibly non-positive form $(dd^{c}\phi)^{n}/n!$,
as is usually shown by combining the Riemann-Roch theorem with Kodaira
vanishing. The main point of the present paper is to show the corresponding
\emph{local} statement. In fact, all three measures will be shown
to coincide with the measure \[
1_{D}(dd^{c}\phi)^{n}/n!\]
where $1_{D}$ is the characteristic function of the set \begin{equation}
D=\{\phi_{e}=\phi\}\subset X\label{eq:D in intro}\end{equation}
 In the case when the metric $\phi$ has a semi-positive curvature
form, $\phi_{e}=\phi,$ i.e. the set $D$ equals all of $X.$

\begin{figure}
\begin{center}\includegraphics{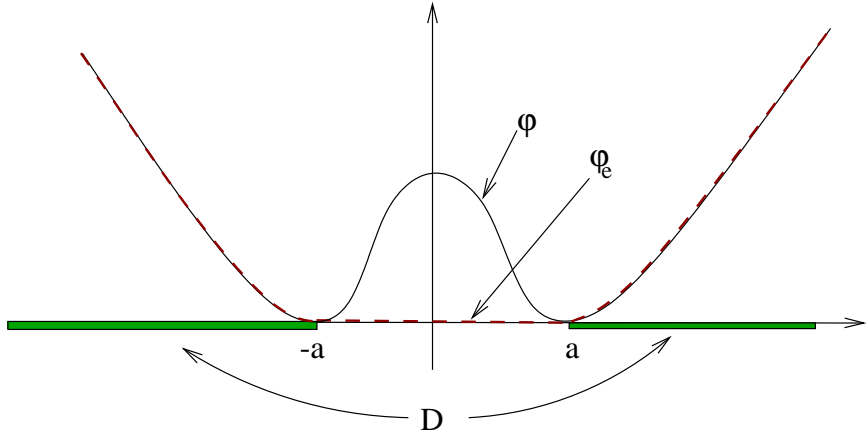}\end{center}

\caption{\label{cap:intro}In the toric case (see example \ref{exa:toric})
the metric $\phi$ is represented by a function on $\R^{n}$ and the
graph of the equilibrium potential $\phi_{e}$ (the dashed line) is
obtained as the \emph{convex hull} of the graph of $\phi.$ The set
$D$ is then the projection of the set where the two graphs coincide. }
\end{figure}

Before turning to the the statement of the main general results, note
that when $L$ is not ample, the main new feature is that the equilibrium
metric $\phi_{e}$ will usually have \emph{singularities} (i.e. points
where it is equal to $-\infty)$ and its curvature $dd^{c}\phi_{e}$
is a positive \emph{current}. However, as is well-known it does give
a metric on $L$ with \emph{minimal} singularities. Such metrics play
a key role in complex geometry (compare remark \ref{rem:minimal}).
Similarly, the Bergman metric $k^{-1}\textrm{ln\,$K_{k}(x,x)$}$ is
singular along the base locus $\textrm{Bs}(\left|\textrm{$kL$}\right|)$
of $L,$ i.e. along the commun zero-locus of the sections in $H^{0}(X,L^{k}).$
Still, the convergence results refered to above in the case when $L$
is ample will be shown to hold provided that the measures are extended
by zero over the singularities. Then integrating over $X$ gives new
proofs of Boucksom's version of Fujita's approximation theorem for
the \emph{volume} of the line bundle $L$ \cite{fu,bo} and its interpretation
in terms of intersection of zero-sets of sections in $H^{0}(X,L^{k})$
by Demailly-Ein-Lazarsfeld \cite{d-e-l}. Finally, the asymptotic
properties of the full Bergman kernel $K_{k}(x,y)$ are studied.

The present approach to the Bergman kernel asymptotics is based on
the use of {}``local holomorphic Morse-inequalities'', which are
local version of the global ones introduced by Demailly \cite{d1}.
These inequalities are then combined with some $L^{2}-$estimates
and global pluripotential theory, the pluripotential part being based
on the recent work \cite{g-z} by Guedj-Zeriahi. Conversely, it turns
out that several basic, but non-trivial, results in pluripotential
theory may obtained as consequences of the Bergman kernel asymptotics
(compare for example remark \ref{rem:vanishing of equil meas by bergm}). 

A crucial step is to first show the $\mathcal{C}^{1,1}$-regularity
of the equilibrium metric $\phi_{e}$ on the complement of the augmented
base-locus $\B_{+}(L),$ which should be of independent interest.

\subsection{Statement of the main results}

Assume that $(L,\phi)$ and $(F,\phi_{F})$ are smooth Hermitian line
bundles over $X$ and denote by $E(k)$ the twisted line bundle $L^{k}\otimes F.$
The equilibrium measure on $X$ (associated to the smooth metric $\phi$
on $L)$ is defined as the positive measure\[
\mu_{\phi}:=1_{U(L)}(dd^{c}\phi_{e})^{n}/n!,\]
where $U(L)$ is the open set in $X$ where $\phi$ is locally bounded
(see section \ref{sec:Equilibrium-measures-for}). The first theorem
to be proved is used to express $\mu_{\phi}$ in terms of $(dd^{c}\phi)^{n}$
on the set $D$ (formula \ref{eq:D in intro} above).

\begin{thm}
\label{thm:reg intro}Suppose that $L$ is a big line bundle and that
the given metric $\phi$ on $L$ is smooth (i.e. in the class $\mathcal{C}^{2}).$
Then the equilibrium metruc $\phi_{e}$ is locally in the class $\mathcal{C}^{1,1}$
on $X-\B_{+}(L)$ i.e. $\phi_{e}$ is differentiable and all of its
first partial derivatives are locally Lipschitz continuous there.
Moreover, the equilibrium measure satisfies \[
\mu_{\phi}n!=1_{X-\B_{+}(L)}(dd^{c}\phi_{e})^{n}=1_{D}(dd^{c}\phi)^{n}=1_{D\cap X(0)}(dd^{c}\phi)^{n}\]
in the sense of measures, where $X(0)$ is the set where $dd^{c}\phi>0$
.
\end{thm}
The regularity theorem is essentially optimal (compare the examples
\ref{exa:toric} and \ref{exa:blow-up}). The next theorem which should
be considered as the main result of this paper gives that, in general,
the measure $k^{-n}B_{k}\omega_{n}$ introduced above (formula \ref{eq:intro B meas})
converges to the equilibrium measure $\mu_{\phi}$. 

\begin{thm}
\label{thm:B in L1 intro}Let $B_{k}$ be the Bergman function of
the Hilbert space $\mathcal{H}(X,E(k)).$ Then\begin{equation}
k^{-n}B_{k}(x)\rightarrow1_{D\cap X(0)}\det(dd^{c}\phi)(x)\label{eq:l1 conv of B}\end{equation}
 for almost any $x$ in $X,$ where $X(0)$ is the set where $dd^{c}\phi>0$
and $D$ is the set \ref{eq:def of D}. Moreover, the following weak
convergence of measures holds:\[
k^{-n}B_{k}\omega_{n}\rightarrow\mu_{\phi},\]
 where $\mu_{\phi}$ is the equilibrium measure. 
\end{thm}
The Bergman function $B_{k}$ may be interpreted as a {}``dimensional
density'' of the Hilbert space $\mathcal{H}(X,E(k)).$ The asymptotic
(normalized) dimension of $\mathcal{H}(X,L^{k})$ is called the \emph{volume}
of a line bundle $L$ \cite{la}: \begin{equation}
\textrm{Vol}(L):=\limsup_{k}k^{-n}\dim H^{0}(X,L^{k})\label{def: vol}\end{equation}
Integrating the convergence of the Bergman kernel in the previous
theoorem now gives the following version of Fujita's approximation
theorem \cite{fu,bo} (compare remark \ref{rem:rel to fujita} for
a comparison with closely related expressions of $\textrm{Vol}(L)).$ 

\begin{cor}
\label{cor:vol as eq intro}The volume of a line bundle $L$ is given
by the total mass of the equilibrium measure: \begin{equation}
\textrm{Vol}(L)=\int_{X}\mu_{\phi}\label{eq:cor vol}\end{equation}
and $\textrm{Vol}(L)=0$ precisely when $L$ is not big.
\end{cor}
The following theorem gives, in particular, the weak convergence on
$X$ of the $k$ th Bergman volume forms (extended by zero over the
base-locus of $E(k)).$ 

\begin{thm}
\label{thm:ln K intro}Let $K_{k}$ be the Bergman kernel of the Hilbert
space $\mathcal{H}(X,E(k).$ Then the following convergence of Bergman
metrics holds: \[
k^{-1}\phi_{k}\rightarrow\phi_{e}\]
uniformly on any fixed compact subset $\Omega$ of $X-\B_{+}(L).$
More precisely, \begin{equation}
e^{-k(\phi-\phi_{e})}C_{\Omega}^{-1}\leq B_{k}\leq C_{\Omega}k^{n}e^{-k(\phi-\phi_{e})}\label{eq:intro Bk exp small intro}\end{equation}
Moreover, the corresponding $k$ th Bergman volume forms converge
to the equilibrium measure: \[
1_{X-\textrm{Bs}(\left|\textrm{$E(k)$}\right|)}(dd^{c}(k^{-1}\textrm{ln\,$K_{k}(x,x)))^{n}/n!\rightarrow\mu_{\phi}$}\]
 weakly as measures on $X.$ 
\end{thm}
The weak convergence in the theorem above on $X-\B_{+}(L)$ is a consequence
of the uniform convergence of the Bergman metrics $k^{-1}\phi_{k}$
on compacts of $X-\B_{+}(L).$ But to get the weak convergence on
all of $X$ theorem \ref{thm:B in L1 intro} (or rather its corollary
\ref{cor:vol as eq intro}) is invoked.

For an ample line bundle $L$ it is a classical fact that the volume
$\textrm{Vol}(L)$ may be expressed as an intersection number $L^{n}.$
More generally, for any line bundle $L$ over $X$ the intersection
of the zero-sets of $n$ {}``generic'' sections in $H^{0}(X,L^{k})$
with $X-\textrm{Bs}(\left|kL\right|)$ (the complement of the commun
zero-locus of all sections) is a finite number of points. The number
of points is called the \emph{moving intersection number} and is denoted
by $(kL)^{[n]}.$ The following corollary was first obtained in \cite{d-e-l}
from Fujita's approximation theorem (see \cite{la} for further references).
The proof given here combines theorem \ref{thm:ln K intro} with properties
of zeroes of {}``random sections'' \cite{sz2}. 

\begin{cor}
If $L$ is a big line bundle then \[
\textrm{Vol}(L)=\lim_{k\rightarrow\infty}\frac{(kL)^{[n]}}{k^{n}}\]

\end{cor}
The final two theorems concern the full Bergman kernel $K_{k}(x,y).$
First, the weak convergence of the squared point-wise norm of the
$K_{k}(x,y)$ is obtained: 

\begin{thm}
\label{thm:k as meas intro}Let $L$ be a line bundle and let $K_{k}$
be the Bergman kernel of the Hilbert space $\mathcal{H}(X,E(k)).$
Then \[
\begin{array}{lr}
k^{-n}\left|K_{k}(x,y)\right|_{k\phi}^{2}\omega_{n}(x)\wedge\omega_{n}(y)\rightarrow\Delta\wedge\mu_{\phi}\end{array},\]
as measures on $X\times X$, in the weak {*}-topology, where $\Delta$
is the current of integration along the diagonal in $X\times X.$
\end{thm}
Then a generalization of the Tian-Zelditch-Catlin expansion \cite{ze}
for a globally positively curved line bundle is shown to hold for
any (big) Hermitian line bundle $L$ over a compact manifold $X:$

\begin{thm}
\label{thm:asymp expansion intro}Let $L$ be a line bundle and let
$K_{k}$ be the Bergman kernel of the Hilbert space $\mathcal{H}(X,E(k)).$
Any interior point in $D\cap X(0)-\B_{+}(L)$ has a neighbourhood
$U$ where $K_{k}(x,y)e^{-k\phi(x)/2}e^{-k\phi(y)/2}$ (with $x,y$
in $U$) admits an asymptotic expansion as\begin{equation}
k^{n}(\det(dd^{c}\phi)(x)+b_{1}(x,y)k^{-1}+b_{2}(x,y)k^{-2}+...)e^{k\phi(x,y)},\label{eq:exp in prop}\end{equation}
where $b_{i}$ are global well-defined functions expressed as polynomials
in the covariant derivatives of $dd^{c}\phi$ (and of the curvature
of the metric $\omega$) which can be obtained by the recursion given
in \cite{b-b-s}.
\end{thm}
Note that \ref{eq:intro Bk exp small intro} in theorem \ref{thm:ln K intro}
implies that $K_{k}(x,y)e^{-k\phi(x)/2}e^{-k\phi(y)/2}=\left|K_{k}(x,y)\right|_{k\phi}$
is exponentially small as soon as $x$ or $y$ is in the complement
of $D.$ 

\begin{rem}
The assumption on $\phi$ may be relaxed to assuming that $\phi$
is in the class $\mathcal{C}^{1,1}.$ For example, the proof of the
regularity theorem \ref{thm:reg intro} still goes through and the
local Morse inequalities (lemma \ref{lem:(Local-Morse-inequalities)})
still apply (almost everwhere on $X).$ Moreover, all results remain
true (with essentially the same proofs) if $X$ is only assumed to
be Moishezon, i.e. bimemorphically equivalent to a projective manifold
(or equivalently, if $L$ carries some big line bundle). However,
in the remaining cases one would have to prove that $1_{D}\det(dd^{c}\phi)=0$
almost everywhere on $X.$ For example, if $(X,\omega)$ is a Kähler
manifold then this would follow from the following conjecture:
\end{rem}
\begin{conjecture}
Let $\omega'$ be a smooth form cohomologous to the Kähler form $\omega.$
Then the global extremal function $V_{X,\omega'}$ associated to $(X,\omega')$
\cite{g-z} (locally expressed as $\phi'_{e}-\phi',$ where $\omega'=dd^{c}\phi')$
is in the class $\mathcal{C}^{1,1}.$ 
\end{conjecture}

\subsection{Further comparison with previous results}

The present paper can be seen as a global geometric version of the
situation recently studied in \cite{berm4}, where the role of the
Hilbert space $\mathcal{H}(X,L^{k})$ was played by the space of all
polynomials in $\C^{n}$ of total degree less than $k,$ equipped
with a weighted norm (compare section \ref{sec:Examples}). For references
to previous works concerning the complex plane see \cite{berm4}.
The proof of the $\mathcal{C}^{1,1}$-regularity of the equilibrium
metric (on the complement of the augmented base locus of $L)$ is
partly modeled on the proof of Bedford-Taylor \cite{b-t,ko} for $\mathcal{C}^{1,1}-$regularity
of the solution of the Dirichlet problem (with smooth boundary data)
for the complex Monge-Ampere equation in the unit-ball in $\C^{n}.$
The result should also be compared to various $\mathcal{C}^{1,1}-$results
for boundary value problems for complex Monge-Ampere equations on
manifolds with boundary \cite{ch,c-t}, intimately related to the
study of the geometry of the space of Kähler metrics on a Kähler manifold
(see also \cite{p-s,bern2} for other relations to Bergman kernels
in the latter context). However, the present situation rather corresponds
to a \emph{free} boundary value problem (compare remark \ref{rem:free}). 

Further references and comments on the relation to the study of random
polynomials (and holomorphic sections), random eigenvalues of normal
matrices and various diffusion-controlled growth processes studied
in the physics literature can be found in \cite{berm4}.

\subsection{Further generalizations}

In a sequel to this paper \cite{berm5} \emph{subspace} and \emph{restricted}
\emph{versions} of the results in this paper will be obtained. The
subspace version is a generalization of the case when the Hilbert
space $\mathcal{H}(X,L^{k})$ is replaced by the subspace of all sections
vanishing to high order along a fixed divisor in $X$ considered in
the preprint \cite{berm4b}. Fixing a singular metric $\phi_{s}$
on $L$ (with analytic singularities) the Hilbert space $\mathcal{H}(X,L^{k})$
is replaced with the subspace of all global holomorphic sections of
the twisted multiplier ideal sheaf $\mathcal{O}(L^{k}\otimes\mathcal{I}(k\phi_{s}))$
(i.e. the space of all sections $f_{k}$ such that the point-wise
norm $\left|f_{k}\right|^{2}e^{-k\phi_{s}}$ is locally integrable)
equipped with the Hilbert subspace norm in $\mathcal{H}(X,L^{k})$
(i.e. the norm induced by the smooth metric $\phi).$ Similarly, the
equilibrium metric $\phi_{e}$ is replaced by the metric obtained
by further demanding that $\widetilde{\phi}\leq\phi_{s}+C$ (i.e.
that $\widetilde{\phi}$ be more singular than $\phi_{s}$) in the
definition \ref{eq:extem metric} of $\phi_{e}.$ As a special case
new proofs of the results of Shiffman-Zelditch \cite{sz5} about Hilbert
spaces of polynomials with coefficients in a scaled Newton polytope
are obtained.

\emph{The restricted versions} are obtained by fixing an $m-$dimension
complex submanifold $V$ of $X$ and replacing $\mathcal{H}(X,L^{k})$
with the restricted space $\mathcal{H}(X,L^{k})_{V}$ equipped with
the {}``restricted norm'' obtained by integrating sections over
$V.$ Similarly, the equilibrium metric $\phi_{e}$ is replaced by
the metric defined on the restricted line bundle $L_{V}$ by only
demanding that $\widetilde{\phi}\leq\phi$ on $V$ in the definition
\ref{eq:extem metric}. The corresponding Bergman kernel asymptotics
can then be seen as local metrized versions of the very recent result
in \cite{e-l-m--} concerning a generalized Fujita approximation theorem
for the restricted volume (i.e. the asymptotic normalized dimension
of $\mathcal{H}(X,L^{k})_{V}).$

\begin{acknowledgement}
It is a pleasure to thank Jean-Pierre Demailly and Sebastian Boucksom
for several illuminating discussions. Also thanks to Frederic Faure
for his interest in this work and for drawing the pictures. This work
was supported by a Marie Curie Intra European Fellowship.
\end{acknowledgement}

\subsection{\label{sub:Notation-and-setup}General notation%
\footnote{general references for this section are the books \cite{gr-ha,de4}.%
}}

Let $(L,\phi)$ be an Hermitian holomorphic line bundle over a compact
complex manifold $X.$ The fixed Hermitian fiber metric on $L$ will
be denoted by $\phi.$ In practice, $\phi$ is considered as a collection
of local \emph{smooth} functions. Namely, let $s^{U}$ be a local
holomorphic trivializing section of $L$ over an open set $U$ then
locally, $\left|s^{U}(z)\right|_{\phi}^{2}=:e^{-\phi^{U}(z)},$ where
$\phi^{U}$ is in the class $\mathcal{C}^{2},$ i.e. it has continuous
derivatives of order two. If $\alpha_{k}$ is a holomorphic section
with values in $L^{k},$ then over $U$ it may be locally written
as $\alpha_{k}=f_{k}^{U}\cdot(s^{U})^{\otimes k},$ where $f_{k}^{U}$
is a local holomorphic function. In order to simplify the notation
we will usually omit the dependence on the set $U.$ The point-wise
norm of $\alpha_{k}$ may then be locally expressed as\begin{equation}
\left|\alpha_{k}\right|_{k\phi}^{2}=\left|f_{k}\right|^{2}e^{-k\phi}.\label{eq:ptwise norm}\end{equation}
 The canonical curvature two-form of $L$ is the global form on $X,$
locally expressed as $\partial\overline{\partial}\phi$ and the normalized
curvature form $i\partial\overline{\partial}\phi/2\pi=dd^{c}\phi$
(where $d^{c}:=i(-\partial+\overline{\partial})/4\pi)$ represents
the first Chern class $c_{1}(L)$ of $L$ in the second real de Rham
cohomology group of $X.$ The curvature form of a smooth metric is
said to be \emph{positive} at the point $x$ if the local Hermitian
matrix $(\frac{\partial^{2}\phi}{\partial z_{i}\partial\bar{z_{j}}})$
is positive definite at the point $x$ (i.e. $dd^{c}\phi_{x}>0).$
This means that the curvature is positive when $\phi(z)$ is strictly
\emph{plurisubharmonic} i.e. strictly subharmonic along local complex
lines. We let \[
X(0):=\left\{ x\in X:\, dd^{c}\phi_{x}>0\right\} \]
 More generally, a metric $\phi'$ on $L$ is called (possibly) \emph{singular}
if $\left|\phi'\right|$ is locally integrable. Then the curvature
is well-defined as a $(1,1)-$current on $X.$ The curvature current
of a singular metric is called \emph{positive} if $\phi'$ may be
locally represented by a plurisubharmonic function (in particular,
$\phi'$ takes values in $[-\infty,\infty[$ and is upper semi-continuous
(u.s.c)). In particular , any section $\alpha_{k}$ as above induces
such a singular metric on $L,$ locally represented by $\phi'=\frac{1}{k}\ln\left|f_{k}\right|^{2}.$
If $Y$ is a complex manifold we will denote by $PSH(Y)$ and $SPSH(Y)$
the space of all plurisubharmonic and strictly plurisubharmonic functions,
respectively. 

Fix another line bundle $F$ with a smooth metric $\phi_{F}$ and
consider the following sequence of Hermitian holomorphic line bundles:
\[
E(k)=(L^{k}\otimes F,k\phi+\phi_{F})\]
 Fixing an Hermitian metric two-form $\omega$ on $X$ (with associated
volume form $\omega_{n})$ the Hilbert space $\mathcal{H}(X,E(k))$
is defined as the space $H^{0}(X,E(k))$ with the norm 

\begin{equation}
\left\Vert \alpha_{k}\right\Vert _{k\phi}^{2}(=\int_{X}\left|f_{k}\right|^{2}e^{-(k\phi(z)+\phi_{F})}\omega_{n}),\label{eq:norm restr}\end{equation}
using a suggestive notation in the last equality (compare formula
\ref{eq:ptwise norm}).

\section{Preliminaries: positivity and base loci}

Let $L$ be a holomorphic line bundle over a compact projective Hermitian
manifold $(X,\omega).$

\subsection{Positivity for line bundles and singular metrics}

The following notions of positivity will be used in the sequel \cite{de5}: 

\begin{defn}
The line bundle $L$ is said to be

\emph{$(i)$ pseudo-effective} if it admits a metric $\phi'$ with
positive curvature current: \[
dd^{c}\phi'\geq0\]
$(ii)$ \emph{big} if it admits a metric $\phi'$ with \emph{strictly
positive curvature current:} \begin{equation}
dd^{c}\phi'\geq\epsilon\omega\label{eq: big metric}\end{equation}

$(iii)$ \emph{ample} if it admits a \emph{smooth} metric $\phi'$
with \emph{strictly positive curvature form:} \[
dd^{c}\phi'_{x}>0\]
for all $x$ in $X.$
\end{defn}

\subsection{Base loci}

For each fixed $k,$ the base locus of the line bundle $E(k)$ is
defined as \[
\textrm{Bs}(\left|\textrm{$E(k)$}\right|)=\bigcap_{f_{k}\in H^{0}(X,E(k))}\{ f_{k}=0\}\]
(or as the corresponding ideal). \emph{The stable base locus} $\B(L)$
of a line bundle $L$ is defined \cite{la} as the following analytic
subvariety of $X:$ \[
\B(L):=\bigcap_{k>0}\textrm{Bs}(\left|kL\right|)=\bigcap_{f_{k}\in H^{0}(X,L^{k}),k\in\N}\{ f_{k}=0\}\]
In other words, a point $x$ is in $\B(L)$ precisely when there is
some section $f_{k}$ in $H^{0}(X,L^{k}),$ for some $k,$ which is
non-vanishing at $x.$ Moreover, the \emph{augumented base locus}
$\B_{+}(L)$ is defined in the following way \cite{la}. Fix an ample
line bundle on $X.$ Then \[
\B_{+}(L):=\B(L-\epsilon A)\]
 for any sufficiently small rational number $\epsilon$ (suitable
interpreted using additive notation for tensor products). We will
have great use for the following equivalent analytic definition of
$\B_{+}(L)$ introduced in \cite{bo2}) (and there called the \emph{non-Kähler
locus}) \begin{equation}
X-\B_{+}(L)=\{ x\in X:\,\,\exists\,\textrm{big\, metric $\phi'$ on $L,\,$smooth\, at\,$x$\} ,}\label{eq: boucksoms augm base}\end{equation}
 in the sense that $\phi'$ satisfies \ref{eq: big metric} and is
smooth on some neighbourhood of the point $x.$%
\footnote{the condition that $\phi'$ be smooth at $x$ may be replaced by haveing
$0$ Lelong number at $x.$ %
} In fact, the equivalence of the definitions is a direct consequence
of theorem \ref{thm:extension} below. It amounts to showing that
$\B_{+}(L)$ is the intersection of all effective divisors $E_{k}(=\{ f_{k}=0\})$
appearing in a {}``Kodaira decomposition'' \begin{equation}
L^{k}=A\otimes[E_{k}],\label{eq:kod decomp}\end{equation}
 for some positive natural number $k.$ The point is that give any
such decomposition \begin{equation}
\phi_{+}:=\frac{1}{k}\ln\left|f_{k}\right|^{2}+\phi_{A}\label{eq:big metric in Kod decomp}\end{equation}
(where $\phi_{A}$ is a fixed smooth metric with positive curvature
on $A)$ is a metric on $L$ with \emph{strictly} positive curvature
current such that $\phi_{+}$ is \emph{smooth on $X-E.$} The reason
for calling $\B_{+}$ the augmented base locus is that\[
\B(L)\subseteq\B_{+}(L),\,\,\,\,\,\bigcap_{k>0}\textrm{Bs}(\left|\textrm{$E(k)$}\right|)\subseteq\B_{+}(L)\]

\begin{rem}
\label{rem:augm base}It is well-known \cite{bo2} that a line bundle
$L$ is ample precisely when $\B_{+}(L)=\emptyset$ and $L$ is big
precisely when $\B_{+}(L)\neq X.$ In particular, $L$ is non-ample
and big pricisely when $\B_{+}(L)$ is a non-empty analytic subvarity
of $X).$ 
\end{rem}

\subsection{Extension of sections}

Let $(L',\phi_{L'})$ and $(F',\phi_{F'})$ be line bundles with singular
metrics such that \begin{equation}
dd^{c}\phi_{L'}\geq\omega/C\,\,\,\textrm{and\,\,\,}dd^{c}\phi_{F'}\geq-C\omega\label{eq:cond on prime metrics}\end{equation}
for some positive number $C.$%
\footnote{$(F',\phi_{F'})$ could also be replaced by a smooth Hermitian holomorphic
vector bundle.%
} Recall the following celebrated theorem about $L^{2}-$estimates
for the $\overline{\partial}-$equation, which is the basic analytical
tool in the theory \cite{de4}.

\begin{thm}
\label{thm:Kodaira-H=ED=AF=BF}(Kodaira-Hörmander-Demailly). Let $\phi_{L'}$
and $\phi_{F'}$ be (singular) metrics as in \ref{eq:cond on prime metrics}.
Take $k'$ sufficiently large. Then for any $\overline{\partial}-$closed
$(0,1)-$form $g$ with values in $L'^{k'}\otimes F'$ such that \[
\left\Vert g\right\Vert _{k'\phi_{L'}+\phi_{F'}}^{2}<\infty,\]
 there is a section $u$ with values in $L'^{k'}\otimes F'$ such
that \[
\overline{\partial}u=g\,\,\,\textrm{and\,\,\,}\left\Vert u\right\Vert _{k'\phi_{L'}+\phi_{F'}}^{2}\leq\frac{C}{k'}\left\Vert g\right\Vert _{k'\phi_{L'}+\phi_{F'}}^{2},\]
where the constant $C$ is independent of $k'$ and $g.$
\end{thm}
The following extension theorem is a well-known asymptotic version
of the Ohsawa-Takegoshi theorem \cite{de5}:

\begin{thm}
\label{thm:extension}Let $\phi_{L'}$ and $\phi_{F'}$ be (singular)
metrics as in \ref{eq:cond on prime metrics}. Fix a point $x$ in
$X$ and take $k'$ sufficiently large. Then any element in $(L^{'k'}\otimes F')_{x}$
such that \[
\left|\alpha_{k'}\right|_{k'\phi_{L'}+\phi_{F'}}^{2}(x)<\infty\]
extends to a section in $H^{0}(X,L'^{k}\otimes F')$ such that \begin{equation}
\left\Vert \alpha_{k'}\right\Vert _{k'\phi_{L'}+\phi_{F'}}^{2}\leq C\left|\alpha_{k'}\right|_{k'\phi_{L'}+\phi_{F'}}^{2}(x)\label{eq:uppsk o-t}\end{equation}
 where the constant $C$ is independent of $k'$ and $\alpha_{k'}(x).$ 
\end{thm}
The previous theorem will be used to extend sections from $X-\B_{+}(L).$
The point is that given any (reasonable) metric $k\phi_{L}+\phi_{F}$
on $L^{k}\otimes F$ the following simple lemma provides a {}``strictly
positively curved perturbation'' $\psi_{k}$ to wich theorem \ref{thm:Kodaira-H=ED=AF=BF}
and \ref{thm:extension} apply:

\begin{lem}
\label{lem:positivization}Let $\phi_{L}$ be a singular metric on
$L$ with positive curvature current and let $\phi_{F}$ be a smooth
metric on the line bundle $F.$ For any given point $x_{0}$ in $X-\B_{+}(L)$
and $k$ sufficiently large there is a singular metric $\psi_{k}$
on $L^{k}\otimes F$ such that $\psi_{k}=k'\phi_{L'}+\phi_{F'}$ with
$\phi_{L'}$ and $\phi_{F'}$ as in formula \ref{eq:cond on prime metrics}
and \begin{equation}
\sup_{U(x_{0})}\left|(k\phi_{L}+\phi_{F})-\psi_{k}\right|\leq C_{x_{0}}\label{eq:sup in lemma pos}\end{equation}
for some neighbourhood $U(x_{0})$ of $x_{0}.$ Moreover, if $\phi_{L}$
is a metric with minimal singularities (see remark \ref{rem:minimal}),
for example the equilibrium metric $\phi_{e}$ (definition \ref{eq:extem metric}),
then it may further be assumed that \begin{equation}
\psi_{k}\leq(k\phi_{L}+\phi_{F})\label{eq:upper in lemma pos}\end{equation}
 on all of $X.$
\end{lem}
\begin{proof}
By the definition of $X-\B_{+}(L')$ there is a metric $\phi_{+}$
on $L,$ smooth in some $U(x_{0})$ with \emph{strictly} positive
curvature current on $X.$ Let $(F',\phi_{F'})=(L^{k-k'}\otimes F,(k-k')\phi_{L}+\phi_{F})$
and $(L',\phi')=(L,\phi_{+})$, for $k'$ fixed. Then $L^{k}\otimes F=L^{k'}\otimes F'$
gets an induced metric\begin{equation}
\psi_{k}=(k-k')\phi_{L}+\phi_{F'}+k'\phi_{+}\label{pf of thm lnk: metric}\end{equation}
of the required form for $k'(\leq k)$ sufficiently large, satisfying
\ref{eq:sup in lemma pos}. Finally, if $\phi_{L}$ is a metric with
minimal singularities we may assume that $\phi_{+}\leq\phi_{L}$ after
substracting a sufficiently large constant from $\phi_{+}.$ Then
\ref{eq:upper in lemma pos} clearly holds.
\end{proof}

\section{\label{sec:Equilibrium-measures-for}Equilibrium measures for line
bundles}

Let $L$ be a line bundle over a compact complex manifold $X.$ Given
a smooth metric $\phi$ on $L$ the corresponding {}``equilibrium
metric'' $\phi_{e}$ is defined as the envelope 

\begin{equation}
\phi_{e}(x)=\sup\left\{ \widetilde{\phi}(x):\,\widetilde{\phi}\in\mathcal{L}_{(X,L)},\,\widetilde{\phi}\leq\phi\,\,\textrm{on$\, X$}\right\} .\label{eq:extem metric}\end{equation}
 where $\mathcal{L}_{(X,L)}$ is the class consisting of all (possibly
singular) metrics on $L$ with positive curvature current. Then $\phi_{e}$
is also in the class $\mathcal{L}_{(X,L)}$ (proposition \ref{pro:of equil}
below). The Monge-Ampere measure $(dd^{c}\phi_{e})^{n}/n!$ is well-defined
on the open set \[
U(L):=\{ x:\,\phi_{e}\,\,\textrm{is\, bounded\, on$\, U(x)\}$},\]
where $U(x)$ is some neighbourhood of $x$ (see \cite{b-t,kl,g-z}
for the definition of the Monge-Ampere measure of a locally bounded
metric or plurisubharmonic function). \emph{The equilibrium measure}
(associated to the metric $\phi)$ is now defined as 

\begin{equation}
\mu_{\phi}:=1_{U(L)}(dd^{c}\phi_{e})^{n}/n!\label{eq:def of mu}\end{equation}
and is hence a positive measure on $X.$ Consider the following set
\begin{equation}
D:=\{\phi_{e}=\phi\}\subset X,\label{eq:def of D}\end{equation}
which is closed by $(i)$ in the following proposition. 

\begin{prop}
\label{pro:max prop of equil}The following holds

$(i)$ $\phi_{e}$ is in the class $\mathcal{L}_{(X,L)}.$

$(ii)$ $1_{U(L)}(dd^{c}\phi_{e})^{n}/n!=0$ on $X-D.$

$(iii)$ $D\subset\{ x:\,$ $dd^{c}\phi_{x}\geq0\}.$ 
\end{prop}
\begin{proof}
$(i)$ is obtained by combining theorem 5.2 (2) and proposition 5.6
in \cite{g-z}. 

The property $(ii)$ is proved precisely as in the local theory in
$\C^{n}$ (compare lemma 2.3 in the appendix of \cite{s-t}). Indeed,
it is enough to prove the vanishing on any small ball in $X-D.$ For
an alternative proof, using Bergman kernels, see the remark below.

To prove $(iii)$ fix a point $x$ where $dd^{c}\phi_{x}<0.$ Then
there is a positive number $\epsilon$ and local coordinates $z$
centered at $x$ such that $(\frac{\partial^{2}\phi}{\partial\zeta\partial\zeta})(\zeta,0,...,0))\leq\epsilon$
(with $z_{1}=\zeta)$ for $\zeta$ in the unit-disc $\Delta.$ Now
take a candidate $\widetilde{\phi}$ for the sup \ref{eq:extem metric}
and let $\psi_{\epsilon}(\zeta):=\widetilde{\phi}(\zeta,0,...,0)-\phi(\zeta,0,...,0)-\epsilon\left|\zeta\right|^{2}.$
Then $\psi_{\epsilon}(\zeta)\leq-\epsilon$ on $\partial\Delta$ and
$\frac{\partial^{2}\psi_{\epsilon}}{\partial z_{1}\partial\bar{z_{1}}}\geq0$
on $\Delta.$ Hence, the submean inequality for subharmonic functions
(or the maximum principle) applied to $\psi_{\epsilon}$ gives $\widetilde{\phi}(x)-\phi(x)=\psi_{\epsilon}(0)\leq-\epsilon.$
Taking the sup over all candidates $\widetilde{\phi}$ then gives
$\phi_{e}(x)-\phi(x)\leq-\epsilon$ which proves the proposition.
\end{proof}
\begin{rem}
\label{rem:vanishing of equil meas by bergm}As will be shown below
$\phi_{e}$ is in the class $\mathcal{C}^{1}$ and its second derivatives
exist almost everywhere and are locally bounded on $X-\B_{+}(L).$
Hence, one could also take $1_{X-\B_{+}(L)}\det(dd^{c}\phi_{e})\omega_{n}$
as a somewhat more concrete definition of the equilibrium measure
on $X$ (which a posteriori anyway gives the same measure on $X,$
according to theorem \ref{thm:reg} below). It is interesting to see
that the vanishing of $1_{X-\B_{+}(L)}\det(dd^{c}\phi_{e})\omega_{n}$
on $X-D$ (corresponding to $(ii)$ in the previous proposition) becomes
a corollary of the proof of theorem \ref{thm:B in L1} (which is independent
of the proof of $(ii)$ in the previous proposition).
\end{rem}
\begin{prop}
\label{pro:of equil}The following properties of equlibrium metrics
hold

$(i)$~ $(m\phi)_{e}=m\phi_{e}$

$(ii)\,$Let $\phi_{A}$ be a metric on a line bundle $A$ such that
$dd^{c}\phi_{A}\geq0.$ Then \[
D_{\phi}\subseteq D_{\phi+\phi_{A}}\]
$(iii)$ Assume that $L$ is big and let $\phi_{F}$ be a smooth metric
on a line bundle $F.$ Then for any compact subset of $X-\ \B_{+}(L)$
there is a constant $C$ such that \[
\frac{C}{m}-\phi_{e}\leq(\phi+\frac{1}{m}\phi_{F})_{e}-\frac{1}{m}\phi_{F}\leq\phi_{e}+\frac{C}{m}\]
for all positive natural numbers $m.$ 
\end{prop}
\begin{proof}
$(i)$ is trivial. For $(ii)$ note that $\phi_{e}+\phi_{A}$ is a
contender for the sup in the definition of $(\phi+\phi_{A})_{e}.$
Hence, for $x$ in $D_{\phi}$ we get\[
\phi(x)+\phi_{A}(x)=\phi_{e}(x)+\phi_{A}(x)\leq(\phi+\phi_{A})_{e}(x).\]
This means that $x$ is also in $D_{\phi+\phi_{A}},$ proving $(ii).$ 

To prove $(iii)$ fix a compact set $\Omega$ in $X-\ \B_{+}(L)$
and a metric $\phi_{+}$ on $L$ with strictly \emph{}positive curvature
current, such that $\phi_{+}$ is smooth $\Omega$ with $\phi_{+}\leq\phi$
on $X.$ Let us first prove one side of the inequality, i.e. \begin{equation}
\phi_{m}:=(\phi+\frac{1}{m}\phi_{F})_{e}-\frac{1}{m}\phi_{F}\leq\phi_{e}+\frac{C}{m}\label{pf of iii lemma proper eq}\end{equation}
 To this end first note that $\phi_{m}$ is a metric on $L$ and there
is clearly a constant $C$ such that \[
\phi_{m}\leq\phi,\,\,\, dd^{c}\phi_{m}\geq-\frac{C}{m}dd^{c}\phi_{+}.\]
 Now let $\phi_{m,+}:=(1-\frac{C}{m})\phi_{m}+\frac{C}{m}\phi_{+},$
defining another metric on $L.$ Note that \[
dd^{c}\phi_{m,+}\geq(1-\frac{C}{m})(-\frac{C}{m}dd^{c}\phi_{+})+\frac{C}{m}dd^{c}\phi_{+}=(\frac{C}{m})^{2}dd^{c}\phi_{+}>0\]
and since also $\phi_{m,+}\leq\phi,$ the extremal definition of $\phi_{e}$
forces $\phi_{m,+}\leq\phi_{e}.$ But since $\phi_{+}$ is smooth
on the compact set $\Omega$ this proves \ref{pf of iii lemma proper eq}.
The other side of the inequality is obtained from \ref{pf of iii lemma proper eq}
applied to $\phi'=\phi+\frac{1}{m}\phi_{F}$ and $\phi'_{F}=-\phi_{F}.$
\end{proof}
The next theorem gives the regularity properties of the equilibrium
metric $\phi_{e}.$ 

\begin{thm}
\label{thm:reg}Suppose that $L$ is a big line bundle and that the
given metric $\phi$ on $L$ is smooth (i.e. in the class $\mathcal{C}^{2}).$
Then 

(a) $\phi_{e}$ is locally in the class $\mathcal{C}^{1,1}$ on $X-\B_{+}(L)$
i.e. $\phi_{e}$ is differentiable and all of its first partial derivatives
are locally Lipschitz continuous there.

(b) The Monge-Ampere measure of $\phi_{e}$ on $X-\B_{+}(L)$ is absolutely
continuous with respect to any given volume form and coincides with
the corresponding $L_{loc}^{\infty}$ $(n,n)-$form obtained by a
point-wise calculation: \begin{equation}
(dd^{c}\phi_{e})^{n}=\det(dd^{c}\phi_{e})\omega_{n}\label{eq:ptwise repr of equil meas}\end{equation}

(c) the following identity holds almost everywhere on the set $D-\B_{+}(L),$
where $D=\{\phi_{e}=\phi\}:$ \begin{equation}
\det(dd^{c}\phi_{e})=\det(dd^{c}\phi)\label{eq:monge on D}\end{equation}
More precisely, it holds for all $x$ in $D-\B_{+}(L)-G,$ where $G$
is the set defined in the proof of $(c).$

(d) Hence, the following identity between measures on $X$ holds:
\begin{equation}
n!\mu_{\phi}=1_{X-\B_{+}(L)}(dd^{c}\phi_{e})^{n}=1_{D}(dd^{c}\phi)^{n}=1_{D\cap X(0)}(dd^{c}\phi)^{n}\label{eq:d in theorem on equil}\end{equation}

\end{thm}
\begin{proof}
$(a),(b)$ and $(c)$ will be proven in the subsequent section. To
prove $(d)$ first observe that the last equality in \ref{eq:d in theorem on equil}
follows immediately from proposition \ref{pro:max prop of equil}$(iii)$.
The second equality in \ref{eq:d in theorem on equil} is obtained
by combining $(c)$ in the stated theorem with the vanishing in proposition
\ref{pro:max prop of equil}$(ii).$ Alternatively, the vanishing
is obtained by combining the bound \ref{eq:strong morse} applied
to $\phi_{e}=\phi'_{e}$ with theorem \ref{thm:B in L1}, giving \[
1_{X-\B_{+}(L)}(dd^{c}\phi_{e})^{n}=1_{D}(dd^{c}\phi_{e})^{n}\]
 Finally, to obtain the vanishing of $(dd^{c}\phi_{e})^{n}$ on $U(L)\cap\B_{+}(L)$
one can use the well-known local fact \cite{kl} that the Monge-Ampere
measure of a locally bounded psh function integrates to zero over
any pluripolar set (in particular over any local piece of $\B_{+}(L)).$ 
\end{proof}
\begin{rem}
\label{rem:minimal}For a general line bundle $L$ the equilibrium
metric $\phi_{e}$ is an example of a metric with \emph{minimal singularities}
in the sense that for any other metric $\phi'$ in $\mathcal{L}_{(X,L)}$
there is a constant $C$ such that \[
\phi'\leq\phi_{e}+C\]
on $X$ (when such an inequality holds $\phi'$ is said to be \emph{more
singular} than $\phi_{e})$ Such metrics play a key role in complex
geometry \cite{de5}.
\end{rem}

\subsection{The proof of $\mathcal{C}^{1,1}-$regularity away from the augmented
base locus}

As in \cite{berm4}, where the manifold $X$ was taken as $\C^{n},$
the proof is modeled on the proof of Bedford-Taylor \cite{b-t,ko,de3}
for $\mathcal{C}^{1,1}-$regularity of the solution of the Dirichlet
problem (with smooth boundary data) for the complex Monge-Ampere equation
in the unit-ball in $\C^{n}.$ However, as opposed to $\C^{n}$ and
the unit-ball a generic compact Kähler manifold $X$ has no global
holomorphic vector fields. In order to circumvent this difficulty
we will reduce the regularity problem on $X$ to a problem on the
manifold $Y,$ where $Y$ is the total space of the dual line bundle
$L^{*},$ identifying the base $X$ with its embedding as the zero-section
in $Y.$ To any given (possibly singular) metric $\phi$ on $L$ we
may associate the logarithm $\chi_{\phi}$ of the {}``squared norm
function'' on $Y,$ where locally \begin{equation}
\chi_{\phi}(z,w)=\ln(\left|w\right|^{2})+\phi(z),\label{eq:h is local coord}\end{equation}
in terms of local coordinates $z_{i}$ on the base $X$ and $w$ along
the fiber of $L^{*}.$ In this way we obtain a bijection \begin{equation}
\mathcal{L}_{(X,L)}\leftrightarrow\mathcal{L}_{Y},\,\,\,\phi\mapsto\chi_{\phi}\label{eq:bijection}\end{equation}
where $\mathcal{L}_{Y}$ is the class of all positively logarithmically
2-homogeneous plurisubharmonic functions on $Y:$ \begin{equation}
\mathcal{L}_{Y}:=\{\chi\in PSH(Y):\,\chi(\lambda\cdot)=\ln(\left|\lambda\right|^{2})+\chi(\cdot)\},\label{eq:class ly}\end{equation}
using the natural multiplicative action of $\C^{*}$ on the fibers
of $Y$ over $X.$ Now we define \begin{equation}
\chi_{e}:=\sup\left\{ \chi\in\mathcal{L}_{Y}:\,\chi\leq\chi_{\phi}\,\,\textrm{on$\, Y$}\right\} ).\label{def: h e}\end{equation}
 Then clearly, $\chi_{e}$ corresponds to the equilibrium metric $\phi_{e}$
under the bijection \ref{eq:bijection}. 

In the following we will denote by $\pi$ the projection from $Y$
onto $X$ and by $j$ the natural embedding of $X$ in $Y.$ We will
fix a point $y_{0}$ in $Y-(j(X)\cup\pi^{-1}(\B_{+}(L)).$ Then there
is a divisor $E$ (appearing in a Kodaira decomposition as in formula
\ref{eq:kod decomp}) such that $y_{0}$ is in $Y-(j(X)\cup\pi^{-1}(E)).$
Since clearly $(k\phi)_{e}=k\phi_{e}$ we may (since we are only interested
in the regularity of $\phi_{e})$ without loss of generality assume
that $k$ appearing in formula \ref{eq:kod decomp} is equal to one.
Moreover, we fix an associated metric $\phi_{+}$ (as in formula \ref{eq:big metric in Kod decomp})
that we write as \[
\phi_{+}:=\phi_{A}+\ln\left|e\right|^{2},\]
 using the suggestive notation $e$ for the defining section of $E,$
and $e(z)$ for any local representative. We may assume that $\phi_{+}\leq\phi_{e}.$
Later we will also assume the normalization $\chi_{\phi_{e}}(y_{0})=0.$

\subsubsection*{Existence of vector fields}

The next lemma provides the vector fields needed in the modification
of the approach of Bedford-Taylor.

\begin{lem}
\label{lem:exist of vf}Assume that the line bundle $L$ is big. For
any given point $y_{0}$ in $Y-[j(X)\cup\pi^{-1}(E)]$ there are global
holomorphic vector fields $V_{1},...V_{n+1}$ (i.e. elements of $H^{0}(Y,TY)$)
such that their restriction to $y_{0}$ span the tangent space $TY_{y_{0}}.$
Moreover, given any positive integer $m$ the vector fields may be
chosen to satisfy 

\begin{equation}
(i)\,\left|V_{i}\right|\leq C_{m}(\left|w\right|)^{m},\,\,\,\,\,(ii)\,\left|V_{i}\right|\leq C_{m}(\left|e(z)\right|)^{m}\label{eq:bound in lemma vf}\end{equation}
locally on the set $\{\chi_{_{\phi_{+}}}\leq1\}$ in $Y$ (in the
following we will fix some $V_{i}$ corresponding to $m=2).$
\end{lem}
\begin{proof}
First note that $Y$ may be compactified by the following fiber-wise
projectivized vector bundle:\[
\widehat{Y}:=\P(L^{*}\oplus\underline{\C}),\]
 where $\underline{\C}$ denotes the trivial line bundle over $X.$
Denote by $\mathcal{O}(1)$ the line bundle over $\widehat{Y}$ whose
restriction to each fiber (i.e. a one-dimensional complex space $\P^{1})$
is the induced hyperplane line bundle. Next, we equip the line bundle
\[
\widehat{L}:=(\pi^{*}(L^{k_{0}})\otimes\mathcal{O}(1))\]
 over $\widehat{Y},$ where $\pi$ denotes the natural projection
from $\widehat{Y}$ to $X,$ with a metric $\widehat{\phi}$ defined
in the following way. First fix a smooth metric $\phi_{E}$ on the
line bundle $[E]$ over $X.$ Then \[
\phi_{+,k_{0}}:=(\phi_{A}-\frac{1}{\sqrt{k_{0}}}\phi_{E})+(1+\frac{1}{\sqrt{k_{0}}})\ln\left|e\right|^{2}\]
is a metric on $L$ such that $dd^{c}\phi_{+}\geq dd^{c}\phi_{A}/2$
for $k_{0}>>1.$ Hence, 

\[
\widehat{\phi}:=\pi^{*}(k_{0}\phi_{+,k_{0}})+\ln(1+e^{\chi_{\phi}})\]
is a metric on $\widehat{L}$ over $Y$ (extending to $\widehat{Y}$)
with strictly positive curvature currrent on $\widehat{Y},$ if $k_{0}>>1.$
Now fix a point $y_{0}$ in $(Y-(j(X)\cup\pi^{-1}(E).$ We can apply
$(ii)$ in theorem \ref{thm:extension} to the bundle $\widehat{L}{}^{k_{1}}\otimes T\widehat{Y}$
over $\widehat{Y}$ for $k_{1}$ sufficiently large. Restricting to
$Y$ in $\widehat{Y}$ then gives that $TY\otimes\pi^{*}(L){}^{k_{0}k_{1}}$
is globally generated on $Y$ (since $\mathcal{O}(1)$ is trivial
on $Y).$ Finally, observe that \[
\pi^{*}(L)=(\pi^{*}(L^{*}))^{-1}=[X]^{-1},\]
 where $[X]$ is the divisor in $Y$ determined by the embedding of
$X$ as the base. Indeed, $X$ is embedded as the zero-set of the
tautological section of $\pi^{*}(L^{*})$ over $Y(=L^{*}).$ Hence,
the sections of $TY\otimes\pi^{*}(L){}^{k_{0}}$ may be identified
with sections in $TY$ vanishing to order $k_{0}$ on $X.$ This proves
$(i)$ in \ref{eq:bound in lemma vf}. Moreover, by construction the
vector fields $V_{i}$ satisfy \[
\left|V_{i}(z,w)\right|\leq C\left(\left|w\right|^{k_{0}}\left|e(z)\right|^{k_{0}(1+\frac{1}{\sqrt{k_{0}}})}\left|w\right|\right)^{k_{1}}\]
on any fixed neighbourhood in $Y$ over the divisor $E$ in $X.$
Choosing $k_{0}$ and $k_{1}$ sufficiently large then gives \[
\left|V_{i}(z,w)\right|\leq C_{m}(\left|w\right|\left|e(z)\right|)^{k_{0}k_{1}+1}\left|e(z)\right|^{m}.\]
Since, the factor $\left|w\right|\left|e(z)\right|$ is bounded on
the set $\{\chi_{\phi_{A}+\ln\left|e\right|^{2}}\leq0\}$ this proves
$(ii)$ in \ref{eq:bound in lemma vf}.
\end{proof}

\subsubsection*{Existence of the flow}

For any given smooth vector field $V$ on $Y$ and compact subset
$K$ of $Y,$ we denote by $\exp(tV)$ the corresponding flow which
is well-defined for any {}``time'' $t$ in $[0,t_{K}],$ i.e. the
family of smooth maps indexed by $t$ such that \begin{equation}
\frac{d}{dt}f(\exp(tV)(y))=df[V]_{\textrm{exp}(tV)(y)}\label{eq:def of flow}\end{equation}
for any smooth function $f$ and point $y$ on $Y.$ We will also
use the notation $\exp(V):=\exp(1V).$ 

Combining the previous lemma with the inverse function theorem gives
local {}``exponential'' holomorphic coordinates centered at $y_{0},$
i.e a local biholomorphism \[
\C^{n+1}\rightarrow U(y),\,\,\,\lambda\mapsto\exp((V(\lambda)(y_{0}),\,\,\, V(\lambda):=\sum\lambda_{i}V_{i})\]
We will write \[
f^{\lambda}=(\exp(V(\lambda))^{*}f\]
for the induced additive action on functions $f$ (where the flow
is defined). Using that the vector fields $V_{i}$ necessarily also
span $TY_{y_{1}}$ for $y_{1}$ close to $y_{0}$ it can be checked
that in order to prove that a function $f$ is locally Lipschitz continuous
on a compact subset of $Y$ it is enough to, for each fixed point
$y_{0},$ prove an estimate of the form \begin{equation}
\left|f^{\lambda}(y_{0})-f(y_{0})\right|\leq C\left|\lambda\right|\label{eq:cond for lip}\end{equation}
for some constant $C$ only depending on the function $f.$ Since
we will later take $f$ to be equal to $\chi_{\phi_{e}}$ we may also,
by homogenity, assume that $\chi_{\phi_{e}}(y_{0})=0.$ In order to
define the flow on a neighbourhood $U$ of the whole levelset $\{\chi_{_{\phi_{e}}}=0\}$
(which is non-compact unless $\phi_{e}$ is locally bounded) we will
use a compactification argument: 

\begin{lem}
\label{lem:exist of flow}There is a positive number $t_{0}$ such
that the flow $\exp(V(\lambda)y)$ exists for any $(\lambda,y)$ such
that $\left|\lambda\right|\leq t_{0}$ and $y$ is in $U:=\{\chi_{_{\phi_{+}}}\leq1\}.$
Moreover, if $\phi'$ is a fixed metric on $L$ such that $\phi'-\delta\ln\left|e\right|^{2}$
is smooth, for some number $\delta,$ then there are constants $C_{\alpha}$
such that\begin{equation}
\left|\partial_{z,w}^{\alpha}(\chi_{\phi'}^{\lambda}-\chi_{\phi'})\right|\leq C_{\alpha}\left|\lambda\right|\label{eq:statement of lemma exist of fl}\end{equation}
 on $U\cap\pi^{-1}(X-E)$ over any fixed $z-$coordinate ball in $X,$
in terms of the real local derivatives of multi order $\alpha$ and
total order less than two. 
\end{lem}
\begin{proof}
Denote by $Y'$ the total space of the ample line bundle $L\otimes[E]^{-1}$
over $X$ and denote by $\pi'$ the corresponding projection onto
$X.$ Then $Y\cap\pi^{-1}(X-E)$ is biholomorphic to $Y'\cap\pi'^{-1}(X-E)$
under the map $\Phi$ which may be locally represented as \begin{equation}
\Phi:(z,w)\mapsto(z',w'):=(z,e(z)w).\label{pf of exist flow: def of isom}\end{equation}
Note that $\Phi$ maps the set $\{\chi_{\phi_{+}}\leq1$ to the set
$\{\chi_{\phi_{A}}\leq1\}.$ Given a vector field $V$ on $Y\cap\pi^{-1}(X-E)$
denote by $V'$ the vector field $\Phi_{*}V$ on $Y'\cap\pi'^{-1}(X-E).$
Now fix a point in $E$ corresponding to $z=0$ in some local coordinates
$(z,w)$ for $Y.$ Then the following local bound holds on the set
$\{\chi_{\phi_{A}}\leq1\}$ in $Y'$ over $X-E:$ \begin{equation}
\left|V_{i}'(z',w')\right|\leq C\left|e(z')\right|.\label{pf: exist of flow local: local bd}\end{equation}
Indeed, writing $V_{i}(z,w)=v_{i,z}(z,w)\frac{\partial}{\partial z}+v_{i,w}(z,w)\frac{\partial}{\partial w}$
and similarly for $V_{i}'(z',w')$ gives \begin{equation}
v_{i,z'}'(z',w')=v_{i,z}(z,w),\,\,\, v_{i,w'}'(z',w')=v_{i,z}(z,w)\frac{\partial e(z)}{\partial z}\frac{1}{e(z)}w'+v_{i,w}(z,w)e(z)\label{eq:pf exist flow}\end{equation}
 Hence, chosing vector fields $V_{i}$ corresponding to $m=2$ in
$(ii)$ in lemma \ref{lem:exist of vf} ensures that \ref{pf: exist of flow local: local bd}
holds. 

Now fix a vector $\lambda$ and write $\lambda=t\sigma,$ where $\left|\lambda\right|=t.$
Let $V$ be the vector field on $Y$ defined by the relation $V(\lambda)=tV.$
By the local bounds \ref{pf: exist of flow local: local bd}, $V'$
extends to a holomorphic vector field on the set $\{\chi_{\phi_{A}}\leq1\}$
in $Y'$ such that $V'$ vanishes identically on $\pi'^{-1}(E).$
Since, $\{ y':\,\chi_{\phi_{A}}(y')\leq1\}$ is compact in $Y'$ the
flow $\textrm{exp}(tV')(y')$ is well-defined for $\left|t\right|\leq t_{0}.$
Moreover, since $V'$ vanishes identically on $\pi'^{-1}(E)$ the
set $\pi'^{-1}(X-E)$ is invariant under the flow. By the isomorphism
$\Phi$ in \ref{pf of exist flow: def of isom} this proves the existence
of the flow stated in the lemma under the assumption that $y\in\pi{}^{-1}(X-E).$
But by the local bound $(ii)$ in lemma \ref{lem:exist of vf} the
flow does extend holomorphically over $\pi{}^{-1}(X-E).$ 

To prove \ref{eq:statement of lemma exist of fl}, note that since
$\chi_{\phi'}$ is smooth over $X-E$ the defining property \ref{eq:def of flow}
of the flow gives\begin{equation}
\chi_{\phi'}^{\lambda}-\chi_{\phi'}=\int_{0}^{1}d\chi_{\phi'}[V(\lambda)]_{\textrm{exp}(tV(\lambda)(y)}dt\label{pf of exist of flow: newton}\end{equation}
Hence, since $V(\lambda):=\left|\lambda\right|(\sum\frac{\lambda_{i}}{\left|\lambda\right|}V_{i}),$
the constant $C_{0}$ in \ref{eq:statement of lemma exist of fl}
may be taken to be \[
C_{0}=\sup_{y\in U\cap\pi^{-1}(X-E)}\left|d\chi_{\phi'}[V_{i}]_{y}\right|<\infty\]
 To see that $C_{0}$ is finite it is, by the compactness of $X,$
enough to prove the bound over any $z-$coordinate ball in $X.$ First
consider local coordinates $(z,w)$ on $Y$ where $z$ is centered
at a point in $X-E.$ Then \begin{equation}
\left|d\chi_{\phi'}[V_{i}]\right|=2\left|\frac{1}{w}dw[V_{i}]+\frac{\partial\phi'(z)}{\partial z}dz[V_{i}]\right|\leq C\label{pf of lem exist of fl: local bd b}\end{equation}
using the bound $(i)$ in lemma \ref{lem:exist of vf} for the first
term in the right hand side above and the assumption that $\phi'$
is smooth on $X-E$ for the second term. Finally, consider the case
when $z$ is centered at a point in $E.$ Then \[
\left|\frac{\partial\phi'(z)}{\partial z}\right|\leq C+\delta\left|\frac{1}{e(z)}\frac{\partial e(z)}{\partial z}\right|\leq C+C'\left|\frac{1}{e(z)}\right|.\]
Hence, the bound \ref{pf of lem exist of fl: local bd b} still holds,
using that $m(:=2)\geq1$ in the local bound $(ii)$ on $V_{i}$ in
lemma \ref{lem:exist of vf}.

Finally note that the bounds for $C_{\alpha}$ when the total degree
of $\alpha$ is two may be obtained in a completely similar maner,
now using that $m=2$ to handle the factors $\left|\frac{1}{e(z)}\right|^{2}$
and $\left|\frac{1}{w}\right|^{2}.$ 
\end{proof}

\subsubsection*{Homogenization}

To a given psh function $g(y)$ (defined on some disc subbundle of
$Y)$ we will associate the following $S^{1}-$invariant psh function:\begin{equation}
\widehat{g}(y):=\textrm{u.s.c}(\sup_{\theta\in[0,2\pi[}g(e^{i\theta}y))\label{def of hat}\end{equation}
using the natural multiplicative action of $\C^{*}$ on the fibers
of $Y$ over $X,$ where u.s.c. denotes the upper-semicontinous regularization
(using that the family $g(e^{i\theta}\cdot)$ of psh functions is
locally bounded from above \cite{kl}). 

The following simple lemma will allow us to {}``homogenize'' in
the normal direction of $Y,$ as well:

\begin{lem}
\label{lem:homo}Suppose that the function $f$ is $S^{1}-$invariant
and psh on some disc subbundle $U$ of $Y$ containing the set $\{ f\leq c\}$
in $(Y-\pi^{-1}(E)),$ where $E$ is an analytic variety in $X.$
Moreover, assume that $(i)$ $f$ is strictly increasing along the
fibers of $Y$ over $X$ and $(ii)$ $f<c$ on the base $X$ and $f>c$
on $\partial U.$ Then there is a function $\widetilde{f}$ in the
class $\mathcal{L}_{Y}$ such that $\widetilde{f}=f$ on the set $\{ f=c\}.$ 
\end{lem}
\begin{proof}
It is enough to construct such a function $\widetilde{f}$ on $Y-\pi^{-1}(E)).$
Indeed, then $\widetilde{f}$ is locally bounded from above close
to $\pi^{-1}(E))$ and hence extends as a unique psh function to $Y$
\cite{kl} (more generally $E$ may be allowed to be locally pluripolar).
Hence we we may without loss of generality assume that $E$ is empty
in the following. First we will show that $f^{-1}(c)$ is an $S^{1}-$subbundle
of $Y$ over $X,$ i.e. the claim that the equation \[
f(\sigma)=c\]
 has a unique solution on $\pi^{-1}(x),$ modulo the action of $S^{1},$
for each fixed point $x$ in $X.$ To this end we identify the restriction
of $f$ to $\pi^{-1}(x)$ with a convex function $g(v)$ of $v=\ln\left|w\right|^{2}$
(using that $f$ is $S^{1}-$invariant and psh). In particular, $g$
is continuous. By the assumption $(i)$ the equation $g(v_{0})=c$
has at most one solution. Moreover, by $(ii)$ and the fact that $g$
is continuous the solution $v_{0}$ does exist. This proves the claim
above. Now define \[
\widetilde{f}(r\sigma):=\ln(r^{2})+f(\sigma).\]
To see that $\widetilde{f}$ is psh, we may, since the problem is
local, assume that $X$ is a ball in $\C_{z}^{n}.$ Write $\widetilde{f}(z,w)=\widetilde{\phi}(z)+\ln\left|w\right|^{2}$
and note that \[
\Omega=U\cap\{\widetilde{\phi}(z)+\ln\left|w\right|^{2}\leq c\}=U\cap\{ f\leq c\}\]
 using that $f$ is strictly increasing along the fibers. Since $f$
is psh it follows that $\partial\Omega$ is pseudoconvex. A classical
result of Bremermann \cite{br} for such Harthogs domains $\Omega$
now implies that $\widetilde{\phi}(z)$ and hence $\widetilde{f}(z,w)$
is psh.%
\footnote{alternatively one can, by local approximation, reduce to the case
when $\widetilde{f}$ is smooth and $c$ regular. Then $dd^{c}\widetilde{f}=(\widetilde{f}/f)dd^{c}f\geq0$
along $T^{1,0}(\partial\Omega)$ and hence, by homogeneity, $\widetilde{f}$
is psh on $Y.$%
}
\end{proof}
\begin{lem}
\label{lem:hom for lamb}For \label{lem:homo for chi}each $\lambda$
(with sufficiently small norm) the function $f=\widehat{\chi_{\phi_{e}}^{\lambda}}$
satisfies the assumptions in the previous lemma with $c=c_{\lambda}=\widehat{\chi_{\phi_{e}}^{\lambda}}(y_{0})$
and $U=\{\chi_{_{\phi_{+}}}\leq1\}.$ 
\end{lem}
\begin{proof}
Given local coordinates on $U,$ let\[
f_{\lambda}(z,w):=\chi_{\phi_{e}}^{\lambda}(z,w)-\ln\left|w\right|^{2}=[\phi_{e}^{\lambda}(z,w)+((\ln\left|w\right|^{2})^{\lambda}-\ln\left|w\right|^{2})].\]
 Since $\phi_{e}^{\lambda}(z,w)$ is psh (a local version of) \ref{eq:statement of lemma exist of fl}
in lemma \ref{lem:exist of flow}) applied to $\phi'=\ln\left|w\right|^{2}$
gives \[
\frac{\partial^{2}f_{\lambda}(z,w)}{\partial w\partial\bar{w}}\geq-C\left|\lambda\right|\]
Note that the constant $C$ may be taken to be independent of the
local coordinates as in the proof of \ref{eq:statement of lemma exist of fl}
in lemma \ref{lem:exist of flow}, since the base $X$ is compact.
Now write 

\[
\widehat{\chi_{\phi_{e}}^{\lambda}}=\widehat{g}-C\left|\lambda\right|\left|w\right|^{2}+\ln\left|w\right|^{2}\]
 with $g=f_{\lambda}+C\left|\lambda\right|\left|w\right|^{2}.$%
\footnote{over a neighbourhood of $E$ we use the local coordinates $(z',w')$
so that the function $\left|w'\right|^{2}$ is uniformly bounded on
$U$. %
} By the maximum principle $\widehat{g}$ in the definition \ref{def of hat}
is always increasing in $v:=\ln\left|w\right|^{2}.$ Hence the (right)
derivative of the convex function $\widehat{\chi_{\phi_{e}}^{\lambda}}$
with respect to $v$ (with the variable $z$ fixed) is positive (and
almost equal to$1)$ for $\lambda$ sufficiently small. This proves
$(i).$

To prove $(ii)$ first observe that \begin{equation}
c_{\lambda}:=\widehat{\chi_{\phi_{e}}^{\lambda}}(y_{0})\leq1/2\label{eq:homo for chi}\end{equation}
 for all $\lambda$ (with sufficiently small norm). Indeed, by upper
semi-continuity $\limsup_{\lambda\rightarrow0}\widehat{\chi_{\phi_{e}}^{\lambda}}(y_{0})\leq\widehat{\chi_{\phi_{e}}}(y_{0})=0.$
Hence, \ref{eq:homo for chi} holds (after possibly replacing the
upper bound $t_{0}$ on $\left|\lambda\right|$ with a smaller number).
Next, note that by the extremal definition of $\phi_{e}$ we have
$\widehat{\chi_{\phi_{e}}^{\lambda}}\geq\chi_{\phi_{+}}^{\lambda}$
on $Y.$ By \ref{eq:homo for chi} and the definition of $U$ it is
hence enough to prove that $\chi_{\phi_{+}}^{\lambda}\geq\chi_{\phi_{+}}-C\left|\lambda\right|$
on $U.$ But this follows from \ref{eq:statement of lemma exist of fl}
in lemma \ref{lem:exist of flow} applied to $\phi'=\phi_{+}.$
\end{proof}

\subsubsection{Proof of $(a)-(c)$ in theorem \ref{thm:reg}}

To prove $(a)$ it is, by the bijection \ref{eq:bijection}, equivalent
to prove that $\chi_{e}$ (defined by \ref{def: h e}) is locally
$\mathcal{C}^{1,1}$ on $Y-[j(X)\cup\pi^{-1}(\B_{+}(L)].$ 

\begin{proof}
\emph{Step1:} \textbf{\emph{}}\emph{$\chi_{e}(:=\chi_{\phi_{e}})$
is locally Lipschitz continuous on $Y-(j(X)\cup\pi^{-1}(\B_{+}(L)).$ }

To see this fix a point $y_{0}$ as above. From the definition \ref{def of hat}
of $\widehat{(\cdot)}$ we have an upper bound \[
\chi_{\phi_{e}}^{\lambda}(y_{0})\leq\widehat{\chi_{\phi_{e}}^{\lambda}}(y_{0})=\widetilde{\widehat{\chi_{\phi_{e}}^{\lambda}}}(y_{0}),\]
 where $\widetilde{\widehat{\chi_{\phi_{e}}^{\lambda}}}$ is the function
in the class $\mathcal{L}_{Y},$ extending $\widehat{\chi_{\phi_{e}}^{\lambda}}$
from the level set \[
M_{\lambda}:=\{ y:\,\widehat{\chi_{\phi_{e}}^{\lambda}}(y)=\widehat{\chi_{\phi_{e}}^{\lambda}}(y_{0})\}\subset U\]
 obtained from lemma \ref{lem:homo for chi}. Since, by definition,
$\chi_{\phi_{e}}^{\lambda}\leq\chi_{\phi}^{\lambda}$ we have the
following bound on the level set $M_{\lambda}:$\begin{equation}
\widetilde{\widehat{\chi_{\phi_{e}}^{\lambda}}}(y)\leq\sup_{\theta\in[0,2\pi]}\chi_{\phi}((\textrm{exp}(V(\lambda))(e^{i\theta}y)\leq\sup_{\theta\in[0,2\pi]}\chi_{\phi}(e^{i\theta}y)+C\left|\lambda\right|,\label{eq:pf thm reg}\end{equation}
using that $\chi_{\phi}$ is smooth in the last inequality so that
\ref{eq:statement of lemma exist of fl} in lemma \ref{lem:exist of flow}
can be applied. Since $\chi_{\phi}$ is $S^{1}-$invariant \ref{eq:pf thm reg}
gives that \begin{equation}
\widetilde{\widehat{\chi_{\phi_{e}}^{\lambda}}}-C\left|\lambda\right|\leq\chi_{\phi}\label{pf thm reg: bound on cand}\end{equation}
 on the level set $M_{\lambda}$ and hence, by homogeneity, on all
of $Y.$ This shows that the function $\widetilde{\widehat{\chi_{\phi_{e}}^{\lambda}}}-C\left|\lambda\right|$
is a contender for the supremum in the definition \ref{def: h e}
of $\chi_{e}$ and hence bounded by $\chi_{e}.$ All in all we get
that \[
\chi_{\phi_{e}}^{\lambda}(y_{0})\leq\widetilde{\widehat{\chi_{\phi_{e}}^{\lambda}}}(y_{0})\leq\chi_{\phi_{e}}(y_{0})+C\left|\lambda\right|.\]
The other side of the inequality \ref{eq:cond for lip} for $f=\chi_{\phi_{e}}^{\lambda}$
is obtained after replacing $\lambda$ by $-\lambda.$ 

\emph{Step2: $d\chi_{e}$ exists and is locally Lipschitz continuous
on $Y-X.$}

Following the exposition in \cite{de3} of the approach of Bedford-Taylor
it is enough to prove the following inequality:

\begin{equation}
(\chi_{\phi_{e}}^{-\lambda}(y_{0})+\chi_{\phi_{e}}^{\lambda}(y_{0}))/2-\chi_{\phi_{e}}(y_{0})\leq C\left|\lambda\right|^{2},\label{eq:pf of thm reg ineq 2}\end{equation}
where the constant only depends on the second derivatives of $\chi_{\phi}.$
Indeed, given this inequality (combined with the fact that $\chi_{\phi_{e}}$
is psh) a Taylor expansion of degree $2$ gives the following bound
close to $y_{0}$ for a local smooth approximation $\chi_{\epsilon}$
of $\chi_{\phi_{e}}:$ \[
\left|D^{2}\chi_{\epsilon}\right|\leq C\]
 where $\chi_{\epsilon}:=\chi_{\phi_{e}}*u_{\epsilon},$ using a a
local regularizing kernel $u_{\epsilon}$ and where $D^{2}\chi_{\epsilon}$
denotes the real local Hessian matrix of $\chi_{\phi_{e}}.$ Letting
$\epsilon$ tend to $0$ then proves Step 2. Finally, to see that
the inequality \ref{eq:pf of thm reg ineq 2} holds we apply the argument
in Step 1 after replacing $\chi_{\phi_{e}}^{\lambda}$ by the psh
function \[
g(y):=(\chi_{\phi_{e}}^{\lambda}(y)+\chi_{\phi_{e}}^{-\lambda}(y))/2\]
 to get the following bound on the level set $\{ y:\,\widehat{g}(y)=\widehat{g}(y_{0})\}:$\[
g(y)\leq\widetilde{\widehat{g}}((y)\leq\sup_{\theta\in[0,2\pi]}(\chi_{\phi}((\textrm{exp}(V(\lambda))(e^{i\theta}y)+\chi_{\phi}(\textrm{exp}(V(-\lambda))(e^{i\theta}y))/2\]
Next, observe that for each fixed $\theta$ the function $\chi_{\phi}(e^{i\theta}y)$
is in the class $\mathcal{C}^{2}.$ Hence, a Taylor expansion of degree
$2$ in the right hand side of formula \ref{pf of exist of flow: newton}
gives \[
\widetilde{\widehat{g}}((y)\leq\sup_{\theta\in[0,2\pi]}((\chi_{\phi}(e^{i\theta}y))+C\left|\lambda\right|^{2})=\chi_{\phi}(y)+C\left|\lambda\right|^{2}\]
 where the constant $C$ may be taken as a constant times $\sup_{U}\left|\left\langle D^{2}\chi_{\phi}V,V\right\rangle _{\eta}\right|$
in terms of some fixed metric $\eta$ on $TY.$ This shows that $\widetilde{\widehat{g}}-C\left|\lambda\right|^{2}$
is a contender for the supremum in the definition \ref{def: h e}
of $\chi_{e}$ and hence bounded by $\chi_{e}.$ All in all we obtain
that \[
g(y_{0})\leq\chi_{\phi}(y_{0})+C\left|\lambda\right|^{2},\]
 which proves the inequality \ref{eq:pf of thm reg ineq 2}, finishing
the proof of Step2.

$(c)$ To see that \ref{eq:monge on D} holds, it is enough to prove
that locally\begin{equation}
\frac{\partial^{2}\phi}{\partial x_{i}\partial x_{j}}(\phi_{e}-\phi)=0\label{pf of thm reg: same deriv}\end{equation}
almost everywhere on $D=\{\phi_{e}=\phi\},$ where $x_{i}$ is a real
coordinate on $\R^{2n}(=\C^{n}).$ To this end let $\psi:=\phi_{e}-\phi$
and let $A:=\{\psi=0,\, d\psi\neq0\}.$ By $(a)$ above $\psi$ is
a $\mathcal{C}^{1}-$function and hence $A$ is a real hypersurface
of codimension $1$ and in particular of measure zero (w.r.t. Lesbegue
meaure). Next, let $f:=d\psi$ (considered as a local map on $\R^{2n})$
and let $B_{1}$ be the set where the derivative $df$ (i.e the matrix
$(df_{1},...,df_{2n})$ does not exist. Since, by $(a)$ $f$ is a
Lipshitz map it is well-known that $B_{1}$ also has measure zero.
Finally, let $B_{2}$ be the set where $f=0,\, df$ exists, but $df\neq0.$
Now using a lemma in \cite{k-s} (page 53) applied to the Lipshitz
map $f$ gives that $B_{2}$ too has measure zero. Finally, let $G:=A\bigcup B_{1}\bigcup B_{2}.$
Then $G$ has measure zero and \ref{pf of thm reg: same deriv} holds
on $(X-\B_{+}(L))-G,$ proving $(c).$ 
\end{proof}
\begin{rem}
\label{rem:free}Suppose that $L$ is ample and fix a smooth metric
$\phi_{+}$ on $L$ with positive curvature. Then $\omega_{+}:=dd^{c}\phi_{+}$
is a Kähler metric on $X$ and the fixed metric $\phi$ on $L$ may
be written as $\phi=u+\phi_{+},$ where $u$ is a smooth \emph{function}
on $X.$ Now the pair $(u_{e},M)$ where $u_{e}:=\phi_{e}-\phi_{+}$
and $M$ is the set $X-D,$ may be interpreted as a {}``weak'' solution
to the following \emph{free} boundary value problem of Monge-Ampere
type%
\footnote{since there is a priori no control on the regularity of the set $M,$
it does not really make sense to write $\textrm{$\partial M$}$ and
the boundary condition should hence be interpreted in a suitable {}``weak''
sense.%
}: \[
\begin{array}{rclr}
(dd^{c}u_{e}+\omega_{+})^{n} & = & 0 & \textrm{on\,$M$}\\
u_{e} & = & u & \textrm{on\,$\partial M$}\\
du_{e} & = & du\end{array}\]
The point is that, since the equations are overdetermined, the set
$M$ is itself part of the solution. In \cite{h-m} the $\mathcal{C}^{1,1}-$regularity
of $\phi_{e}$ in the case when $X=\C$ (corresponding to the setup
in \cite{s-t}) was deduced from the regularity of a free boundary
value problem. 
\end{rem}

\section{\label{sec:Bergman-kernels}Bergman kernel asymptotics}

Recall that $\mathcal{H}(X,E(k))$ denotes the Hilbert space obtained
by equipping the vector space $H^{0}(X,E(k))$ with the norm \ref{eq:norm restr}.
Let $(\psi_{i})$ be an orthonormal base for $\mathcal{H}(X,E(k)).$
The \emph{Bergman kernel} of the Hilbert space $\mathcal{H}(X,E(k))$
is the integral kernel of the orthogonal projection from the space
of all smooth sections with values in $E(k)$ onto $\mathcal{H}(X,E(k)).$
It may be represented by the holomorphic section \begin{equation}
K_{k}(x,y)=\sum_{i}\psi_{i}(x)\otimes\overline{\psi_{i}(y)}.\label{def: K}\end{equation}
 of the pulled back line bundle $E(k)\boxtimes\overline{E(k)}$ over
$X\times\overline{X}.$ The restriction of $K_{k}$ to the diagonal
is a section of $E(k)\otimes\overline{E(k)}$ and we let $B_{k}(x)=\left|K_{k}(x,x)\right|_{k\phi+\phi_{F}}(=\left|K_{k}(x,x)\right|e^{-(k\phi(x)+\phi_{F})})$
be its point wise norm: \begin{equation}
B_{k}(x)=\sum_{i}\left|\psi_{i}(x)\right|_{k\phi+\phi_{F}}^{2}.\label{eq:def of B}\end{equation}
We will refer to $B_{k}(x)$ as the \emph{Bergman function} of $\mathcal{H}(X,E(k)).$
It has the following extremal property:\begin{equation}
B_{k}(x)=\sup\left\{ \left|\alpha_{k}(x)\right|_{k\phi+\phi_{F}}^{2}:\,\,\alpha_{k}\in\mathcal{H}(X,E(k)),\,\left\Vert \alpha_{k}\right\Vert _{k\phi+\phi_{F}}^{2}\leq1\right\} \label{(I)extremal prop of B}\end{equation}
 Moreover, integrating \ref{eq:def of B} shows that $B_{k}$ is a
{}``dimensional density'' of the space $\mathcal{H}(X,L^{k}):$
\begin{equation}
\int_{X}B_{k}\omega_{n}=\dim\mathcal{H}(X,E(k))\label{eq:dim formel for B}\end{equation}
The following {}``local Morse inequality'' estimates $B_{k}$ point-wise
from above for a general bundle:

\begin{lem}
\emph{\label{lem:(Local-Morse-inequalities)}(Local Morse inequalities)}
On any compact complex manifold the following upper bound holds: \[
k^{-n}B_{k}\leq C_{k}1_{X(0)}\det(dd^{c}\phi),\]
 where the sequence $C_{k}$ of positive numbers tends to one and
$X(0)$ is the set where $dd^{c}\phi>0.$
\end{lem}
See \cite{berm1} for the more general corresponding result for $\overline{\partial}-$harmonic
$(0,q)$- forms with values in a high power of an Hermitian line bundle.
The present case (i.e. $q=0)$ is a simple consequence of the mean-value
property of holomorphic functions applied to a poly-disc $\Delta_{k}$
of radius $\ln k/\sqrt{k}$ centered at the origin in $\C^{n}$ (see
the proof in \cite{berm2}). In fact, the proof gives the following
stronger local statement: \begin{equation}
\limsup_{k}k^{-n}\left|f_{k}(z/\sqrt{k})\right|^{2}e^{-(k\phi+\phi_{F})((z/\sqrt{k})}/\left\Vert f_{k}\right\Vert _{k\phi+\phi_{F},\Delta_{k}}^{2}\leq1_{X(0)}(0)\det(dd^{c}\phi),\label{eq:local morse on ball}\end{equation}
 where $f_{k}$ is holomorphic function defined in a fixed neighbourhood
of the origin in $\C^{n}.$

The local estimate in the previous lemma can be considerably sharpened
on the complement of the globally defined set $D$ (formula \ref{eq:def of D}),
as shown by the following lemma:

\begin{lem}
\label{lem:exponentia decay} Let $\phi$ be a smooth metric on a
holomorphic line bundle $L$ over a compact manifold $X.$ 

(i) if $B_{k}$ denotes the Bergman function of the Hilbert space
$\mathcal{H}(X,L^{k}),$ then the following inequality holds on all
of $X$: \begin{equation}
B_{k}k^{-n}\leq C_{k}e^{-k(\phi-\phi_{e})}\label{eq:exp decay}\end{equation}
 where the sequence $C_{k}$ of positive numbers tends to $\sup_{X}\det(dd^{c}\phi).$ 

(ii) If $L$ is big and $B_{k}$ now denotes the Bergman function
of the Hilbert space $\mathcal{H}(X,E(k)),$ then the inequality \ref{eq:exp decay}
holds on any given compact subset of $X-\B_{+}(L)$ if $C_{k}$ is
replaced by a large constant $C$ (depending on the compact set). 

In particular, in both cases\begin{equation}
\lim\int_{D^{c}}k^{-n}B_{k}\omega_{n}=0\label{pf of thm B: int B on D compl}\end{equation}

\end{lem}
\begin{proof}
Let us first prove $(i).$ By the extremal property \ref{(I)extremal prop of B}
of $B_{k}$ it is enough to prove the lemma with $B_{k}k^{-n}$ replaced
by $\left|\alpha_{k}\right|_{k\phi}^{2},$ locally represented by
$\left|f_{k}\right|e^{-k\phi},$ for any element $\alpha_{k}$ in
$\mathcal{H}(X,L^{k})$ with global norm equal to $k^{-n}.$ The Morse
inequalities in the previous lemma give that\[
\left|f_{k}\right|^{2}e^{-k\phi}\leq C_{k}\]
 with $C_{k}$ as in the statement of the present lemma. Equivalently,
\[
\frac{1}{k}\ln\left|f_{k}\right|^{2}-\frac{1}{k}C_{k}\leq\phi\]
 Hence, the singular metric on $L$ determined by $\frac{1}{k}\ln\left|f_{k}\right|^{2}-\frac{1}{k}C_{k}$
is a candidate for the sup in the definition \ref{eq:extem metric}
of $\phi_{e}$ and is hence bounded by $\phi_{e}.$ Thus, \[
B_{k}k^{-n}=\left|f_{k}\right|^{2}e^{-k\phi}\leq C_{k}e^{k\phi_{e}}e^{-k\phi}.\]
The proof of $(ii)$ is completely analogous if one takes into account
that $\ln\left|f_{k}\right|^{2}$ is now a metric on $E(k)=L^{k}\otimes F$
and uses $(ii)$ in proposition \ref{pro:of equil}. Finally, the
vanishing \ref{pf of thm B: int B on D compl} follows from the dominated
convergence theorem (using that the sequence $B_{k}k^{-n}$ is, by
lemma \ref{lem:(Local-Morse-inequalities)}, uniformly bounded on
$X),$ since the right hand side in the previous inequality tends
point-wise to zero precisely on the complement of $D.$
\end{proof}
\begin{rem}
When $L$ is ample the supremum, in the definition of $C_{k},$ in
the previous lemma can be taken over the support of $(dd^{c}\phi_{e})^{n},$
using a max/comparison principle. 
\end{rem}
The following lemma yields a \emph{lower} bound on the Bergman function.

\begin{lem}
\label{lem:lower bound on B}Let $L$ be a big line bundle, then the
following lower bound holds at almost any point $x$ in $D\cap X(0):$
\begin{equation}
\liminf_{k}k^{-n}B_{k}\omega_{n}\geq(dd^{c}\phi)^{n}/n!\label{eq:lower bd on B in lemma}\end{equation}

\end{lem}
\begin{proof}
\emph{Step1: construction of a smooth extremal $\sigma_{k}.$} Fix
a point $x_{0}$ in $D\cap X(0)-\B_{+}(L)-G,$ where $G$ is the set
of measure zero appearing in the proof of $(c)$ in the regularity
theorem \ref{thm:reg}. First note that there is a \emph{smooth} section
$\sigma_{k}$ with values in $E(k)$ such that \begin{equation}
(i)\lim_{k\rightarrow\infty}\frac{\left|\sigma_{k}\right|_{k\phi}^{2}(x_{0})}{k^{n}\left\Vert \sigma_{k}\right\Vert _{k\phi+\phi_{F}}^{2}}\omega_{n}=(dd^{c}\phi)_{x_{0}}^{n},\,\,\,(ii)\left\Vert \overline{\partial}\sigma_{k}\right\Vert _{k\phi_{e}+\phi_{F}}^{2}\leq Ce^{-k/C}\label{pf of lemma lower B: prop of sigma}\end{equation}
To see this take trivializations of $L$ and $F$ and local holomorphic
coordinates $z_{i}$ centered at $x_{0}$ (and orthonormal at $x_{0})$
such that $\phi_{F}(0)=0$ and \begin{equation}
\begin{array}{ccc}
\phi(z)=(\sum_{i=1}^{n}\lambda_{i}\left|z_{i}\right|^{2}+O(\left|z\right|^{3}))\end{array}\label{pf of lemma lower bd B: phi as model}\end{equation}
with $\lambda_{i}$ the positive eigenvalues of $(dd^{c}\phi)_{x_{0}}$
w.r.t the metric $\omega$ \cite{gr-ha}. Fix a smooth function $\chi$
which is constant when $\left|z\right|\leq\delta/2$ and supported
where $\left|z\right|\leq\delta;$ the number $\delta$ will be assumed
to be sufficiently small later on. Now $\sigma_{k}$ is simply obtained
as the local section with values in $L^{k}$ represented by the function
$\chi$ close to $x_{0}$ and extended by zero to all of $X.$ To
see that $(i)$ holds note that, using \ref{pf of lemma lower bd B: phi as model},
\[
\lim_{k\rightarrow\infty}\frac{\left|\sigma_{k}\right|_{k\phi}^{2}(x_{0})}{k^{n}\left\Vert \sigma_{k}\right\Vert _{k\phi+\phi_{F}}^{2}}=\lim_{k\rightarrow\infty}\frac{\chi(0)}{k^{n}\int_{\left|z\right|\leq k^{-1/2}\ln k}e^{-k\sum_{i=1}^{n}\lambda_{i}\left|z_{i}\right|^{2}}\chi(0)\omega_{n}(0)},\]
where $\omega_{n}(0)$ is the Euclidian volume form in $\C^{n}$ (since
$z_{i}$ are assumed to be orthonormal w.r.t $\omega$ at $0).$ Evaluating
the latter Gaussian integral then gives the limit $(1/\pi)^{n}\lambda_{1}\lambda_{2}\cdot\cdot\cdot\lambda_{n},$
proving $(i)$ in \ref{pf of lemma lower B: prop of sigma}. To prove
$(ii)$ in \ref{pf of lemma lower B: prop of sigma}, first note that
\[
\left\Vert \overline{\partial}\sigma_{k}\right\Vert _{k\phi_{e}+\phi_{F}}^{2}\leq C\int_{\delta/2\leq\left|z\right|\leq\delta}e^{-k(\phi(z)+(\phi_{e}(z)-\phi(z))}\omega_{n}(0),\]
 as follows from the definition of $\chi.$ Hence, $(ii)$ follows
from the fact that \begin{equation}
\left|z\right|\leq\delta\Rightarrow\phi(z)+(\phi_{e}(z)-\phi(z))\geq\inf_{i}\lambda_{i}\left|z\right|/2\label{pf of lemma lower bd B: delta}\end{equation}
 for $\delta$ sufficiently small, if $x_{0}$ is in the set $D\cap X(0)-\B_{+}(L)-G.$
To prove \ref{pf of lemma lower bd B: delta} fix $z\neq0$ and let
$\psi(y):=\phi_{e}(y\frac{z}{\left|z\right|})-\phi(y\frac{z}{\left|z\right|}),$
where $y$ is a non-negative number. By the regularity theorem \ref{thm:reg}
$u$ is in the class $\mathcal{C}^{1}.$ Hence, since $x_{0}$ is
in $D$ (where $\psi=0)$ \[
\left|\psi(y)\right|\leq\int_{0}^{y}\left|\frac{d\psi}{ds}\right|(s)ds.\]
Moreover, since $x_{0}$ is in $D-\B_{+}(L)-G,$ we have $\frac{d\psi}{ds}(0)=0$
and $\lim_{s\rightarrow0}\frac{d\psi}{ds}(s)/s=0.$ In particular,
for any given positive number $\epsilon$ we have that \[
\left|\psi(y)\right|\leq\int_{0}^{y}\epsilon sds=\epsilon y^{2}/2\]
 for $y\leq\delta,$ if $\delta$ is sufficiently small. Combining
the latter estimate with \ref{pf of lemma lower bd B: phi as model}
then proves \ref{pf of lemma lower bd B: delta} and finishes the
proof of $(ii)$ in \ref{pf of lemma lower B: prop of sigma}. 

\emph{Step2: perturbation of $\sigma_{k}$ to a holomorphic extremal
$\alpha_{k}.$} Equip $E(k)$ with a {}``strictly positively curved
modification'' $\psi_{k}$ of the metric $k\phi_{e}+\phi_{F}$ furnished
by lemma \ref{lem:positivization}. Let $g_{k}=\overline{\partial}\sigma_{k}$
and let $\alpha_{k}$ be the following holomorphic section \[
\alpha_{k}:=\sigma_{k}-u_{k},\]
 where $u_{k}$ is the solution of the $\overline{\partial}$-equation
in the Hörmander-Kodaira theorem \ref{thm:Kodaira-H=ED=AF=BF} with
$g_{k}=\overline{\partial}\sigma_{k}.$ Hence, \[
\left\Vert u_{k}\right\Vert _{\psi_{k}}\leq C\left\Vert g_{k}\right\Vert _{\psi_{k}}\]
 Next, applying \ref{eq:sup in lemma pos} to the right hand side
above (using that $g_{k}$ is supported on a small neighboorhood of
$x_{0}\in X-\B_{+}(L)$ and then \ref{eq:upper in lemma pos} to the
left hand side above gives \[
\left\Vert u_{k}\right\Vert _{k\phi_{e}+\phi_{F}}\leq C\left\Vert g_{k}\right\Vert _{k\phi_{e}+\phi_{F}}\]
Finally, using that $\phi_{e}\leq\phi$ on all of $X$ in the left
hand side above gives \begin{equation}
\left\Vert u_{k}\right\Vert _{k\phi+\phi_{F}}\leq C\left\Vert g_{k}\right\Vert _{k\phi_{e}+\phi_{F}}\label{eq:uk}\end{equation}
and then $(ii)$ in \ref{pf of lemma lower B: prop of sigma} in the
right hand side gives \[
(a)\,\left\Vert u_{k}\right\Vert _{k\phi+\phi_{F}}\leq Ce^{-k/C},\,\,\,(b)\,\left|u_{k}\right|_{k\phi+\phi_{F}}^{2}(x)\leq C'k^{n}e^{-k/C'},\]
where $(b)$ is a consequence of $(a)$ and the local holomorphic
Morse inequalities \ref{eq:local morse on ball} applied to $u_{k}$
at $z=0.$ Combining $(a)$ and $(b)$ with $(i)$ in \ref{pf of lemma lower B: prop of sigma}
then proves that $(i)$ in \ref{pf of lemma lower B: prop of sigma}
holds with $\sigma_{k}$ replaced by the holomorphic section $\alpha_{k}.$
By the definition of $B_{k}$ this finishes the proof of the lemma.
\end{proof}
\begin{rem}
\label{rem:strong Morse}For \emph{any} line bundle $L$ over $X$
and $\phi'$ a given (singular) metric on $L$ with positive curvature
form Boucksom showed \cite{bo} that 

\begin{equation}
\liminf_{k}k^{-n}\dim H^{0}(X,E(k))\geq\int_{X}((dd^{c}\phi')_{ac})^{n}/n!,\label{eq:strong morse}\end{equation}
 where $(dd^{c}\phi')_{ac}$ denotes the absolutely continious part
of the the current $dd^{c}\phi'.$ Boucksom used Bonavero's strong
Morse inequalities for singular metrics with analytic singularities.
However, when $L$ is big the lower bound \ref{eq:strong morse} follows
from a variant of the proof of the previous lemma. To see this one
first approximates $\phi'$ with a sequence $\phi_{\epsilon}'$ with
\emph{analytic} singularities (as in \cite{bo}) and replace $\phi$
with the metric \[
\phi'_{\epsilon,+}:=\phi'_{\epsilon}(1-\epsilon)+\epsilon\phi_{+}\]
with \emph{strictly} positive curvature form (since $\phi_{+}$ is
of the form \ref{eq: big metric}) and the Hilbert space $\mathcal{H}(X,E(k))$
with the Hilbert space $\mathcal{H}(X,E(k),\phi'_{\epsilon,+})$ whose
norm is defined with respect to the norm induced by the singular metric
$\phi'_{\epsilon,+}.$ %
\footnote{as a vector space the Hilbert space $\mathcal{H}(X,E(k),\phi'_{\epsilon,+})$
consists of sections which have finite norm with respect to the singular
metric $\phi_{\epsilon}'(1-\epsilon)+\epsilon\phi_{+}.$%
}A variant of the proof of the previous lemma then gives a lower bound
on the corresponding Bergman function with $\phi$ replaced by $\phi_{\epsilon}',$
for any point $x$ in the complement of the singularity locus of $\phi_{\epsilon}'.$
Finally, integrating over $X$ and then letting $\epsilon$ tend to
zero gives \ref{eq:strong morse}.
\end{rem}
\begin{thm}
\label{thm:B in L1}Let $B_{k}$ be the Bergman function of the Hilbert
space $\mathcal{H}(X,E(k)).$ Then\begin{equation}
k^{-n}B_{k}(x)\rightarrow1_{D\cap X(0)}\det(dd^{c}\phi)(x)\label{eq:l1 conv of B}\end{equation}
 for almost any $x$ in $X,$ where $X(0)$ is the set where $dd^{c}\phi>0$
and $D$ is the set \ref{eq:def of D}. Moreover, the following weak
convergence of measures holds:\[
k^{-n}B_{k}\omega_{n}\rightarrow\mu_{\phi},\]
 where $\mu_{\phi}$ is the equilibrium measure.
\end{thm}
\begin{proof}
\textbf{Case1:} $L$ is big

First observe that, by the exponential decay in lemma \ref{lem:exponentia decay},
\[
\lim_{k\rightarrow\infty}k^{-n}B_{k}(x)=0,\,\, x\in D^{c}\]
Next, the local Morse inequalities (lemma \ref{lem:(Local-Morse-inequalities)})
give the upper bound in \ref{eq:l1 conv of B} on $k^{-n}B_{k}(x)$
for any $x$ in $D$ and lemma \ref{lem:lower bound on B} gives,
since $L$ is assumed to be big, the lower bound for almost any $x$
in $D,$ finishing the proof of \ref{eq:l1 conv of B}. Finally, the
weak convergence follows from \ref{eq:l1 conv of B} combined with
the uniform upper bound on $k^{-n}B_{k}(x)$ in lemma \ref{lem:(Local-Morse-inequalities)})
(using the dominated convergence theorem on $X).$

\textbf{Case2:} $L$ is pseudo-effective, but not big

First note that the dimension of $H^{0}(X,E(k))$ is of the order
$o(k^{n}),$ i.e. \begin{equation}
\lim_{k\rightarrow\infty}\int_{X}k^{-n}B_{k}\omega_{n}=0.\label{pf of conv of B: vanish}\end{equation}
 To see this one argues by contradiction (compare proposition 6.6$f$
in \cite{de5}) : if the dimension would be of the order $k^{n}$
a standard argument gives that $H^{0}(X,E(k)\otimes A^{-1})$ has
a section for $k$ sufficiently large, if $A$ is any fixed ample
line bundle. But then \[
L^{k}=(A\otimes F^{-1})\otimes E,\]
 where $E$ is an effective divisor. Since we may choose $A$ so that
$A\otimes F^{-1}$ is ample this means that $L$ has a metric with
\emph{strictly} positive curvature form (compare formula \ref{eq:big metric in Kod decomp}),
giving a contradiction. 

Now given \ref{pf of conv of B: vanish}, Fatou's lemma forces \[
\lim_{k\rightarrow\infty}k^{-n}B_{k}\omega_{n}(x)=0\]
for almost all $x$ in $X.$ Finally, to prove the vanishing of the
positive measure $1_{D\cap X(0)}\det(dd^{c}\phi)(x)$ a.e. on $X,$
equip the ample line bundle $A$ with a metric $\phi_{A}$ with positive
curvature form. Then proposition \ref{pro:of equil} $(i)$ gives
that \[
\int_{D_{p\phi}\cap X(0)}(dd^{c}p\phi)^{n}/n!\leq\int_{D_{p\phi+\phi_{A}}\cap X(0)}(dd^{c}(p\phi+\phi_{A})^{n}/n!.\]
 Using $(ii)$ in proposition \ref{pro:of equil} and that $dd^{c}p\phi\leq dd^{c}(p\phi+\phi_{A})$
then gives \begin{equation}
\int_{D\cap X(0)}(dd^{c}\phi)^{n}/n!\leq p^{-n}\int_{D_{p\phi+\phi_{A}}\cap X(0)}(dd^{c}(p\phi+\phi_{A})^{n}/n!=p^{-n}\textrm{Vol}(pL+A),\label{pf of thm B not big}\end{equation}
 where we have applied case $1$ to the big line bundle $pL+A$ in
the last step and used the definition \ref{def: vol} below of the
volume $\textrm{Vol}(L').$ Now fix $\epsilon>0.$ By the {}``continuity''
of the volume function (compare remark \ref{rem:cont}) the right
hand side is bounded by \[
p^{-n}\textrm{Vol}(pL)+\epsilon=\epsilon.\]
Finally, letting $p$ tend to infinity in \ref{pf of thm B not big}
gives that $\int_{D\cap X(0)}(dd^{c}\phi)^{n}/n!$ must vanish, since
$\epsilon$ was arbitrary.

\textbf{Case3:} $L$ is not pseudo-effective (and hence not big)

In this case it follows directly from the definition that $\phi_{e}\equiv-\infty$
and hence the set $D$ is empty.
\end{proof}
\begin{rem}
As shown in the proof of lemma \ref{lem:lower bound on B} the set
$G$ of measure zero where the point-wise convergence in the previous
theorem may fail, can be expressed in terms of the {}``derivatives''
(up to order two) of $\phi_{e}-\phi.$ The proof of lemma \ref{lem:lower bound on B}
also shows that the $L^{1}-$convergence for the corresponding Bergman
functions for weighted polynomials obtained in \cite{berm4} can be
improved to convergence on the complement of $G.$ 
\end{rem}
The volume of a line bundle is defined by the following formula \cite{la}:\begin{equation}
\textrm{Vol}(L):=\limsup_{k}k^{-n}\dim H^{0}(X,L^{k})\label{def: vol}\end{equation}
Integrating the convergence of the Bergman kernel in \ref{thm:B in L1}
combined with the previous corollary gives the following version of
Fujita's approximation theorem \cite{fu,bo}:

\begin{cor}
\label{cor:vol as eq}The volume of a big line bundle $L$ is given
by \begin{equation}
\textrm{Vol}(L)=\int_{X-\B_{+}(L)}(dd^{c}\phi_{e})^{n}/n!\label{eq:cor vol}\end{equation}
and $\textrm{Vol}(L)=0$ precisely when $L$ is not big.
\end{cor}
\begin{rem}
\label{rem:rel to fujita}Assume that $L$ is big. Fujita's approximation
theorem (as formulated in \cite{la,d-e-l}) may be stated as \begin{equation}
\textrm{Vol}(L)=\sup_{A}A^{n}\label{eq:fujitas version}\end{equation}
 where the supremum is taken over the top intersection numbers of
all ample line bundles $A$ occuring in a decompostion \ref{eq:kod decomp}
on some modification of $X.$ In fact, Fujita's proved th\emph{e upper}
bound on $\textrm{Vol}(L)$ and the lower bound is considered to be
substantially easier. The following two analytical versions are due
to Boucksom \cite{bo}: \begin{equation}
\textrm{Vol}(L)=\sup_{\phi'\in\mathcal{L}_{(X,L),a}}\int_{X-\textrm{$\{\phi'=-\infty\}$ }}(dd^{c}\phi')^{n}/n!=\sup_{\phi'\in\mathcal{L}_{(X,L)}}\int_{X}((dd^{c}\phi')_{ac})^{n}/n!\label{eq:bouck versions}\end{equation}
 where $\mathcal{L}_{(X,L),a}$ denotes the subspace of metrics with
analytical singularities and $((dd^{c}\phi')_{ac}$ denotes the absolutely
continious part of the the current $dd^{c}\phi'.$ The equivalence
between \ref{eq:fujitas version} and the first equality in \ref{eq:bouck versions}
is simply obtained by taking a log resolution of the pair $(X,\{\phi'=-\infty\})$
(compare \cite{bo}). The equivalence between \ref{eq:bouck versions}
and corollary \ref{cor:vol as eq} is essentially contained in the
proof of theorem \ref{thm:B in L1}: the lower bound on $\textrm{Vol}(L)$
follows from \ref{eq:strong morse} and the upper bound from the fact
that $(dd^{c}\phi_{e})^{n}$ realizes the supremum in the last equality
in \ref{eq:bouck versions} (and is approximated by $\phi_{j}$ in
$\mathcal{L}_{(X,L),a}$).
\end{rem}
Finally, a remark concerning the continuity of the volume function.

\begin{rem}
\label{rem:cont}As is well-known the volume is continous in the sense
that for any line bundles $L$ and $F$ \[
\lim_{m\rightarrow\infty}m^{-n}\textrm{Vol}(mL+F)=\textrm{Vol}(L).\]
The continuity on the (open) cone of big line bundles is, for example,
a simple consequence of the formula for $\textrm{Vol}$ in corollary
\ref{cor:vol as eq} (compare \cite{bo}). To get continuity up to
the boundary of the big cone (for example that the limit is zero when
$L$ is non-big and pseudo-effective and $F$ is ample) one can replace
lemma \ref{lem:lower bound on B} with the bound \ref{eq:strong morse}
(as in \cite{bo}). For a more direct algebro-geometric argument see
\cite{la} (I prop. 2.2.35). 
\end{rem}

\subsection{The Bergman metric}

The Hilbert space $\mathcal{H}(X,E(k)$ induces a metric on the line
bundle $L$ which may be expressed as \[
\phi_{k}(x):=k^{-1}\textrm{ln\,$K_{k}(x,x)-k^{-1}\phi_{F}(x)$}\]
 When $E(k)=L^{k}$ and $\phi_{F}=0$ this metric is in the class
$\mathcal{L}_{(X,L)}$ and is often referred to as the $k$th Bergman
metric on $L.$ If $L$ is an ample line bundle, then this is the
smooth metric on $L$ obtained as the pull-back of the Fubini-Study
metric on the hyperplane line bundle $\mathcal{O}(1)$ over $\P^{N}(=\P\mathcal{H}(X,L^{k}))$
(compare example \ref{exa:Pn} in section \ref{sec:Examples}) under
the Kodaira map \[
X\rightarrow\P\mathcal{H}(X,L^{k}),\,\,\, y\mapsto(\textrm{$\Psi_{1}(x):\Psi_{2}(x)...:\Psi_{N}(x))$ , }\]
for $k$ sufficiently large, where $\textrm{$(\Psi_{i})$}$ is an
orthonormal base for $\mathcal{H}(X,L^{k})$ \cite{gr-ha}. Note that
the metric $\phi_{k}$ on $L$ is singular precisely on the base locus
 $\textrm{Bs}(\left|\textrm{$E(k)$}\right|)$ and its curvature current
is {}``almost positive'' when $k$ is large. The (almost positive)
measures \[
1_{X-\textrm{Bs}(\left|\textrm{$E(k)$}\right|)}(dd^{c}(k^{-1}\textrm{ln\,$K_{k}(x,x)))^{n}$}/n!\]
 on $X$ will be referred to as the $k$ th Bergman volume forms.

Now we can prove the following theorem:

\begin{thm}
\label{thm:ln K}Let $K_{k}$ be the Bergman kernel of the Hilbert
space $\mathcal{H}(X,E(k).$ Then the following convergence holds:\begin{equation}
k^{-1}\phi_{k}\rightarrow\phi_{e}\label{eq:conv of ln K in theorem ln K}\end{equation}
uniformly on any fixed compact subset $\Omega$ of $X-\B_{+}(L).$
More precisely, \[
e^{-k(\phi-\phi_{e})}C_{\Omega}^{-1}\leq B_{k}\leq C_{\Omega}k^{n}e^{-k(\phi-\phi_{e})}\]
Moreover, the corresponding $k$ th Bergman volume forms converge
to the equilibrium measure: \begin{equation}
1_{X-\textrm{Bs}(\left|\textrm{$E(k)$}\right|)}(dd^{c}(k^{-1}\textrm{ln\,$K_{k}(x,x)))^{n}/n!\rightarrow\mu_{\phi}$}\label{eq:cong of monge in thm ln k}\end{equation}
weakly as measures on $X.$ 
\end{thm}
\begin{proof}
In the following proof it will be convenient to let $C$ denote a
sufficiently large constant (which may hence vary from line to line).
First observe that when $(F,\phi_{F})$ is trivial taking the logarithm
of the inequality \ref{eq:exp decay} in lemma \ref{lem:exponentia decay}
immediately gives the upper bound \[
k^{-1}\textrm{ln\,$K_{k}(x,x)\leq\phi_{e}(x)+C\ln k/k$}\]
and the general case is completely analogous.

To get a lower bound, fix a point $x_{0}$ in $X-\B_{+}(L).$ By the
extremal property \ref{(I)extremal prop of B} it is enough to find
a section $\alpha_{k}$ in $\mathcal{H}(X,E(k)$ such that \begin{equation}
\left|\alpha_{k}(x_{0})\right|_{k\phi_{e}+\phi_{F},\,}\geq1/C\,\,\left\Vert \alpha_{k}\right\Vert _{X,k\phi+\phi_{F}}\leq C.\label{pf of theorem ln K}\end{equation}
 To this end take a {}``strictly positively curved modification''
$\psi_{k}$ of the metric $k\phi_{e}+\phi_{F}$ furnished by lemma
\ref{lem:positivization}. Then the extension theorem \ref{thm:extension}
gives a section $\alpha_{k}$ in $\mathcal{H}(X,E(k)$ such that\[
\left|\alpha_{k}(x_{0})\right|_{\psi_{k}}\geq1/C,\,\,\,\left\Vert \alpha_{k}\right\Vert _{X,\psi_{k}}\leq C,\]
 Applying \ref{eq:sup in lemma pos} to the first inequality above
and then \ref{eq:upper in lemma pos} to the second one proves \ref{pf of theorem ln K}
(also using that by definition $\phi_{e}\leq\phi).$

To prove \ref{eq:cong of monge in thm ln k} first observe that the
weak Monge-Ampere convergence \ref{eq:cong of monge in thm ln k}
on the open set $X-\B_{+}(L)$ follows from the uniform convergence
\ref{eq:conv of ln K in theorem ln K} (see \cite{g-z}). Finally,
by general integration theory it is (by theorem \ref{thm:reg} $(d)$
enough to prove \[
\lim_{k\rightarrow\infty}\int_{X-\textrm{Bs}(\left|\textrm{$E(k)$}\right|)}(dd^{c}(k^{-1}\textrm{ln\,$K_{k}(x,x)))^{n}/n!=\int_{X-\B_{+}(L)}(dd^{c}\phi_{e})^{n}n/!$}\]
To this end first note that the weak Monge-Ampere convergence on $X-\B_{+}(L)$
implies \[
\liminf_{k}\int_{X-\textrm{Bs}(\left|\textrm{$E(k)$}\right|)}(dd^{c}(k^{-1}\textrm{ln\,$K_{k}(x,x)))^{n}/n!\geq\int_{X-\B_{+}(L)}(dd^{c}\phi_{e})^{n}n/!$}\]
 Next, applying formula \ref{eq:bouck versions} to $\phi'=\textrm{ln\,$K_{k}(x,x)$}$
shows that the left hand side (but with $\limsup$ instead) is bounded
by $\textrm{Vol}(L)$ (compare remark \ref{rem:rel to fujita}). By
corollary \ref{eq:cor vol} this proves \ref{eq:cong of monge in thm ln k}
and finishes the proof of the theorem.
\end{proof}
For any line bundle $L$ over $X$ the intersection of the zero-sets
of $n$ {}``generic'' sections in $H^{0}(X,L^{k})$ with $X-\textrm{Bs}(\left|kL\right|)$
is a finite number of points (as follows form Bertini's theorem \cite{gr-ha}).
The number of points is called the \emph{moving intersection number}
and is denoted by $(kL)^{[n]}.$ The following corollary was obtained
in \cite{d-e-l} from Fujita's theorem (see \cite{la} for further
references). 

\begin{cor}
If $L$ is a big line bundle then \[
\textrm{Vol}(L)=\lim_{k\rightarrow\infty}\frac{(kL)^{[n]}}{k^{n}}\]

\end{cor}
\begin{proof}
The proof follows immediately from $(ii)$ in theorem \ref{thm:ln K}
and the following fact:\[
n!(kL)^{[n]}=\int_{X-\B_{k}}(dd^{c}(\textrm{ln\,$K_{k}(x,x)))^{n}.$ }\]
The formula my be deduced by taking a log resolution of the pair $(X,\B_{k}(L))$
(compare \cite{bo}). But it also follows from properties of {}``random
zeroes'', once one accepts that \[
n!(kL)^{[n]}=\int_{X-\B_{k}}Z_{f_{1}}\wedge...\wedge Z_{f_{n}},\]
 where $Z_{f_{1}}$ denotes the integration current determined by
the zero-set of $f_{i},$ is independent of a {}``generic'' tuple
$(f_{1},...f_{n})$ in $(H^{0}(X,L^{k})^{n}$. Indeed, by proposition
2.2 in \cite{s-z 6} the right hand side may be written as \[
\int_{X-\B_{k}}\E(Z_{f_{1}}\wedge...\wedge Z_{f_{n}}),\]
where $\E(Z_{f_{1}}\wedge...\wedge Z_{f_{n}})$ denotes the expectation
value (taking values in the space of measures)%
\footnote{$\E(.)$ denotes integration with respect to the Gaussian probability
measure on the product $(\mathcal{H}(X,L^{k})^{n}$ of Hilbert spaces.%
} of the intersection of the zero currents of $n$ random independent
sections in $\mathcal{H}(X,L^{k}).$ Changing the order of integration
gives\[
\int_{X-\B_{k}}\E(Z_{f_{1}}\wedge...)=\E\int_{X-\B_{k}}(Z_{f_{1}}\wedge...)=n!\E(kL)^{[n]}=n!(kL)^{[n]}\]

\end{proof}

\subsection{The full Bergman kernel }

Combining the convergence in theorem \ref{thm:B in L1} with the local
inequalities \ref{eq:local morse on ball}, gives the following convergence
for the point-wise norm of the full Bergman kernel $K_{k}(x,y).$
The proof is completely analogous to the proof of theorem 2.4 in part
1 of \cite{berm2}. 

\begin{thm}
\label{thm:k as meas}Let $L$ be a line bundle and let $K_{k}$ be
the Bergman kernel of the Hilbert space $\mathcal{H}(X,E(k)).$ Then
\[
\begin{array}{lr}
k^{-n}\left|K_{k}(x,y)\right|_{k\phi}^{2}\omega_{n}(x)\wedge\omega_{n}(y)\rightarrow\Delta\wedge\mu_{\phi}\end{array},\]
as measures on $X\times X$, in the weak {*}-topology, where $\Delta$
is the current of integration along the diagonal in $X\times X.$
\end{thm}
Finally, we will show that around any interior point of the set $D\cap X(0)-\B_{+}(L)$
the Bergman kernel $K_{k}(x,y)$ admits a complete local asymptotic
expansion in powers of $k,$ such that the coefficients of the corresponding
symbol expansion coincide with the Tian-Zelditch-Catlin expansion
(concerning the case when the curvature form of $\phi$ is positive
on all of $X;$ see \cite{b-b-s} and the references therein for the
precise meaning of the asymptotic expansion). We will use the notation
$\phi(x,y)$ for a fixed almost holomorphic-anti-holomorphic extension
of a local representation of the metric $\phi$ from the diagonal
$\Delta$ in $\C^{n}\times\C^{n},$ i.e. an extension such that the
anti-holomorphic derivatives in $x$ and the holomorphic derivatives
in $y$ vanish to infinite order along $\Delta.$ 

\begin{thm}
\label{thm:asymp expansion}Let $L$ be a line bundle and let $K_{k}$
be the Bergman kernel of the Hilbert space $\mathcal{H}(X,L^{k}).$
Any interior point in $D\cap X(0)-\B_{+}(L)$ has a neighbourhood
where $K_{k}(x,y)e^{-k\phi(x)/2}e^{-k\phi(y)/2}$ admits an asymptotic
expansion as\begin{equation}
k^{n}(\det(dd^{c}\phi)(x)+b_{1}(x,y)k^{-1}+b_{2}(x,y)k^{-2}+...)e^{k\phi(x,y)},\label{eq:exp in prop}\end{equation}
where $b_{i}$ are global well-defined functions expressed as polynomials
in the covariant derivatives of $dd^{c}\phi$ (and of the curvature
of the metric $\omega$) which can be obtained by the recursion given
in \cite{b-b-s}.
\end{thm}
\begin{proof}
The proof is obtained by adapting the construction in \cite{b-b-s},
concerning globally \emph{positive} Hermitian line bundles, to the
present situation. The approach in \cite{b-b-s} is to first construct
a {}``local asymptotic Bergman kernel'' with the asymptotic expansion
\ref{eq:exp in prop} close to any point where $\phi$ is smooth and
$dd^{c}\phi>0$. Hence, the local construction applies to the interiour
of the set $D\cap X(0)-\B_{+}(L)$ as well. Then the local kernel
is shown to differ from the true kernel by a term of order $O(k^{-\infty}),$
by solving a $\overline{\partial}$-equation with a good $L^{2}-$estimate.
This is possible since $dd^{c}\phi>1/C$ \emph{globally} in that case.
In the present situation we are done if we can solve \begin{equation}
\overline{\partial}u_{k}=g_{k},\label{eq:inhom dbar}\end{equation}
 where $g_{k}$ is a $\overline{\partial}-$closed $(0,1)-$form with
values in $L^{k},$ supported on the interior of the bounded set $D\cap X(0)-\B_{+}(L)$
with an estimate \begin{equation}
\left\Vert u_{k}\right\Vert _{k\phi+\phi_{F}}\leq C\left\Vert g_{k}\right\Vert _{k\phi+\phi_{F}}\label{eq:horm est2}\end{equation}
 To this end note that proceeding precisely as in step 2 in the proof
of lemma \ref{lem:lower bound on B}, gives according to formula \ref{eq:uk}
a solution $u_{k}$ satisfying \ref{eq:horm est2}, but with $\phi$
replaced with $\phi_{e}$ in the norm of $g_{k}.$ However, since
$\phi_{e}=\phi$ on the set where $g_{k}$ is assumed to be supported
this does prove \ref{eq:horm est2} and hence finishes the proof of
the theorem.
\end{proof}

\section{\label{sec:Examples}Examples}

Finally, we illustrate some of the previous results with the following
examples, which can be seen as variants of the setting considered
in \cite{berm4}.

\begin{example}
\label{exa:Pn}Let $X$ be the $n-$dimensional projective space $\P^{n}$
and let $L$ be the hyperplane line bundle $\mathcal{O}(1)$. Then
$H^{0}(X,L^{k})$ is the space of homogeneous polynomials in of degree
$k$ in the $n+1$ homogeneous coordinates $Z_{0},Z_{1},..Z_{n}.$
The Fubini-Study metric $\phi_{FS}$ on $\mathcal{O}(1)$ may be suggestively
written as $\phi_{FS}(Z)=\ln(\left|Z\right|^{2})$ and the Fubini-Study
metric $\omega_{FS}$ on $\P^{n}$ is the normalized curvature form
$dd^{c}\phi_{FS}.$ Hence the induced norm on $H^{0}(X,L^{k})$ is
invariant under the standard action of $SU(n+1)$ on $\P^{n}.$ We
may identify $\C^{n}$ with the {}``affine piece'' $\P^{n}-H_{\infty}$
where $H_{\infty}$ is the {}``hyperplane at infinity'' in $\C^{n}$
(defined as the set where $Z_{0}=0).$ In terms of the standard trivialization
of $\mathcal{O}(1)$ over $\C^{n}$ (obtained by setting $Z_{0}=1)$
the space $H^{0}(Y,L^{k})$ may be identified with the space of polynomials
$f_{k}(\zeta)$ in $\C_{\zeta}^{n}$ of total degree at most $k$
and the metric $\phi_{FS}$ on $\mathcal{O}(1)$ may be represented
by the function\[
\phi_{FS}(\zeta)=\ln(1+\left|\zeta\right|^{2}).\]
Moreover, any smooth metric on $\mathcal{O}(1)$ may be represented
by a function $\phi(\zeta)$ satisfying the following necessary growth
condition%
\footnote{in order that $\phi$ extend over the hyperplane at infinity to a
\emph{smooth} metric further conditions are needed.%
} \begin{equation}
-C+\ln(1+\left|\zeta\right|^{2})\leq\phi(\zeta)\leq\ln(1+\left|\zeta\right|^{2})+C,\label{eq:growth cond}\end{equation}
 which makes sure that the norm \ref{eq:norm restr}, expressed as
\[
\left\Vert f_{k}\right\Vert _{k\phi}^{2}:=\int_{\C^{n}}\left|f_{k}(\zeta)\right|^{2}e^{-k\phi(\zeta)}\omega_{FS}^{n}/n!\]
 is finite when $f_{k}$ corresponds to a section of the $k$ th power
of $\mathcal{O}(m)$, for $m=1.$ In particular, any smooth compactly
supported function $\chi(\zeta)$ determines a \emph{smooth} perturbation
\begin{equation}
\phi_{\chi}(\zeta):=\phi_{FS}(\zeta)+\chi(\zeta)\label{eq:pert metric}\end{equation}
of $\phi_{FS}$ on $\mathcal{O}(1)$ over $\P^{n},$ to which the
results in section \ref{sec:Equilibrium-measures-for} and \ref{sec:Bergman-kernels}
apply. 
\end{example}
The next class of example is offered by toric varieties.

\begin{example}
\label{exa:toric}Let $\Delta$ be a Delzant polytope in $\R^{n}$
obtained as the convex hull of points in $\Z^{n}$ (see \cite{a})
It induces a triple $(X_{\ \Delta},L_{\ \Delta},\phi_{\Delta}),$
where $X_{\ \Delta},$ is an $n-$dimensional complex compact projective
manifold on which the complex torus $\C^{*n}$ acts effectively with
an open dense orbit and $(L_{\ \Delta},\phi_{\Delta})$ is an Hermitian
positive line bundle, invariant under the action of $T^{n}$ (the
real torus in $\C^{*n}$). The curvature form $dd^{c}\phi$ defines
a $T^{n}-$invariant Kähler metric on $X_{\Delta}.$ Identifying $\C_{z}^{*n}$
with an open dense set in $X_{\ \Delta},$ the space $H^{0}(X_{\ \Delta},L_{\Delta}^{k})$
may be identified with the space spanned by all monomials $z^{\alpha}$
with $\alpha$ a multi-index in the scaled polytope $k\Delta$ and
the metric $\phi_{\Delta}$ may be identified with a plurisubharmonic
function on $\C_{z}^{*n}:$\[
\phi_{\Delta}(z)=\ln(\sum_{\alpha\in\Delta\cap\Z^{n}}\left|z^{\alpha}\right|^{2}).\]
Writing $v_{i}:=\ln(\left|z_{i}\right|^{2})$ identifies $\phi_{\Delta}(z)$
with a convex function on $\R^{n}$ that we, by a slight abusive of
notation, denote by $\phi_{\Delta}(v).$ Now any real smooth compactly
supported function $\chi(v)$ on $\R^{n}$ induces a perturbation
$\phi:=\phi_{\Delta}+\chi,$ yielding a new smooth $T^{n}-$invariant
metric $\phi$ on $L_{\Delta}.$ In this notation the almost everywhere
convergence of the Bergman function $B_{k}$ (theorem \ref{thm:B in L1})
may be written as \[
\sum_{p\in\Delta\cap(\frac{1}{k}\Z)^{n}}\frac{e^{k(\left\langle p,v\right\rangle -\phi(v))}}{\int_{v\in\R^{n}}e^{k(\left\langle p,v\right\rangle -\phi(v))}\det(\frac{\partial^{2}\phi_{\Delta}}{\partial^{2}v})dv}\rightarrow1_{D}(v)\left(\frac{2}{\pi}\right)^{n}\frac{\det(\frac{\partial^{2}\phi}{\partial^{2}v})(v)}{\det(\frac{\partial^{2}\phi_{\Delta}}{\partial^{2}v})(v)}\]
where $\left\langle p,v\right\rangle $ and $dv$ denote the Euclidian
scalar product and volume form, respectively. 

Note that in the logarithmic coordinates $v_{i}$ we have that $(\frac{\partial^{2}\phi}{\partial z_{i}\partial\bar{z_{j}}})=(\frac{\partial^{2}\phi}{\partial v_{i}\partial v_{j}})/4.$
Hence, it follows (more or less from the definition) that the graph
of the equilibrium metric $\phi_{e}$ determined by $\phi$ is simply
the \emph{convex hull} of the graph of $\phi$ considered as a function
of $v.$ In particular, in the {}``generic'' toric case $\phi_{e}$
will not be in the class $\mathcal{C}^{2}.$ Indeed, consider for
example the case when $n=1,$ so that $X_{\Delta}=\P^{1},$ and take
$\phi(v)$ to be an even function with two non-degenerate minima at
$\pm a$ (figure \ref{cap:intro} in section \ref{sec:Introduction}).
Then $\frac{\partial^{2}\phi}{\partial^{2}v}(a)>0,$ but $\frac{\partial^{2}\phi_{e}}{\partial^{2}v}(a-\epsilon)=0$
if $0<\epsilon<2a.$ Note that if we deforme $\phi$ in such a way
that $\phi_{e}$ is unchanged (figure \ref{cap:deformed}), then the
leading Bergman kernel asymptotics for the Bergman function $B_{k}$
given by theorem \ref{thm:B in L1} are also unchanged. %
\begin{figure}
\begin{center}\includegraphics{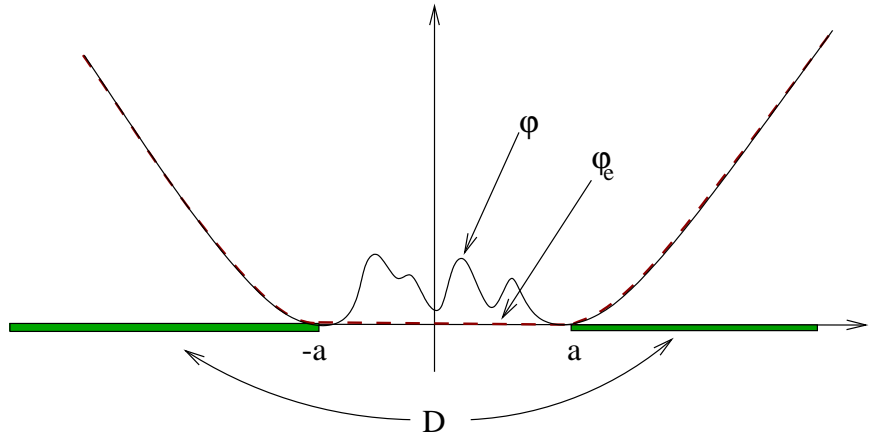}\end{center}

\caption{\label{cap:deformed}The picture represents the graph of a metric
$\phi$ on $\mathcal{O}(2)$ over $\P^{1},$ in logarithmic radial
coordinates, which is a compact perturbation $(\phi=\phi_{\Delta}+\chi)$
of the Fubini-Study metric $\phi_{\Delta}(v)(=\ln(e^{-v}+e^{v})+C).$
Deforming $\phi$ without changing the {}``convex hull'' $\phi_{e}$
does not change the leading Bergman function asymptotics. }
\end{figure}
However, the exponential decay of $B_{k}$ in the complement of the
set $D$ does change. For example, by theorem \ref{thm:ln K}:\[
\lim_{k\rightarrow\infty}\inf_{X}(B_{k})^{1/k}=e^{-\sup_{X}(\phi_{e}-\phi)},\]
i.e. minus the exponential of the limit equals the supremum of the
distance between points on the graph of $\phi$ and points on the
graph of $\phi_{e}.$ In \cite{g-z} the limit in the left hand side
above is called the \emph{Tchebishev constant} and the right hand
side the \emph{Alexander capacity} (associated to the tripple $(X,L,\phi)).$ 
\end{example}
The following basic example of a big (non-ample) line bundle shows
that $\phi_{e}$ may be singular on all of $\B_{+}(L):$

\begin{example}
\label{exa:blow-up}Let $X$ be the blow-up of $\P^{2}$ and denote
by $\pi$ the projection (blow-down map) from $X$ to $\P^{2}.$ Let
$L=\pi^{*}\mathcal{O}(1)\otimes[E],$ where $E$ is the exeptional
divisor. Since, $\int_{E}$$c_{1}(L)<0$ any element of \emph{$\mathcal{L}_{(X,L)}$}
is identically equal to $-\infty$ on $E.$ Moreover, as is well-known
$E=\B_{+}(L).$ In particular, $\phi_{e}\equiv-\infty$ on $\B_{+}(L).$
\end{example}

\end{document}